\documentclass[11pt,twoside,reqno]{amsart}
\usepackage{amsmath,amsthm,amsopn,amssymb,a4wide,times}
\usepackage[mathscr]{eucal} 
\numberwithin{equation}{section}

\newtheorem{theorem}{Theorem}[section]
\newtheorem{lemma}[theorem]{Lemma}
\newtheorem{corollary}[theorem]{Corollary}
\newtheorem{proposition}[theorem]{Proposition}

\newtheorem{definition}[theorem]{Definition}
\newtheoremstyle{note}
     {}
     {}
     {}
     {}
     {\bfseries} 
     {.}
     {.5em}
     {}

\theoremstyle{note}

\newtheorem{remark}[theorem]{Remark}

\newcommand{\ip}[2]{{\left<#1,#2\right>}}
\newcommand{\ran}{{\text{ran}\,}}

\begin{document}
\setcounter{page}{1} \title{The failure of rational dilation on a
  triply connected domain}

\author[M.~A.~Dritschel and S.~McCullough]{Michael A.  Dritschel$^1$
  and Scott McCullough$^2$}

\address{School of Mathematics and Statistics\\
  Merz Court,\\ University of Newcastle upon Tyne\\
  Newcastle upon Tyne\\
  NE1 7RU\\
  UK}

\email{M.A.Dritschel@newcastle.ac.uk}

\address{Department of Mathematics\\
  University of Florida\\
  Box 118105\\
  Gainesville, FL 32611-8105\\
  USA}

\email{sam@math.ufl.edu}

\subjclass[2000]{47A25 (Primary), 30C40, 30E05, 30F10, 46E22, 47A20,
  47A48 (Secondary)}

\keywords{dilations, spectral sets, multiply connected domains, inner
  functions, Herglotz representations, Fay reproducing kernels,
  Riemann surfaces, theta functions, transfer functions,
  Nevanlinna-Pick interpolation}

\thanks{${}^1$Research supported by the EPSRC.  \quad ${}^2$Research
  supported by the NSF}

\date{\today}

\begin{abstract}
  For $R$ a bounded triply connected domain with boundary consisting
  of disjoint Jordan loops there exists an operator $T$ on a complex
  Hilbert space $\mathcal H$ so that the closure of $R$ is a spectral
  set for $T$, but $T$ does not dilate to a normal operator with
  spectrum in $B$, the boundary of $R$.  There is considerable overlap
  with the construction of an example on such a domain recently
  obtained by Agler, Harland and Rafael \cite{AHR} using numerical
  computations and work of Agler and Harland \cite{AH}.
\end{abstract}

\maketitle


\section{Introduction}

Let $R$ denote a domain in $\mathbb C$ with boundary $B$.  Let $X$
denote the closure of $R$.

An operator $T$ on a complex Hilbert space $\mathcal H$ has $X$ as a
spectral set if $\sigma(T)\subset X$ and
\begin{equation*}
  \|f(T)\| \le \|f\|_R = \sup\{|f(z)|: z\in R\}
\end{equation*}
for every rational function $f$ with poles off $X$.  The expression
$f(T)$ may be interpreted either in terms of the Riesz functional
calculus, or simply by writing the rational function $f$ as $pq^{-1}$
for polynomials $p$ and $q$, in which case $f(T)=p(T)q(T)^{-1}$.

The operator $T$ has a normal $B$-dilation if there exists a Hilbert
space $\mathcal K$ containing $\mathcal H$ and a normal operator $N$
on $\mathcal K$ so that
\begin{equation}
  \label{eq:dilation}
  f(T) = P_\mathcal H f(N)|_\mathcal H,
\end{equation}
for every rational function $f$ with poles off the closure of $R$.
Here $P_\mathcal H$ is the orthogonal projection of $\mathcal K$ onto
$\mathcal H$.

It is evident that if $T$ has a normal $B$-dilation, then $T$ has $X$
as a spectral set.  The Sz.-Nagy Dilation Theorem implies that the
closed unit disk is a spectral set for $T$ if and only if $T$ has a
unitary dilation.  A deep result of Agler says that if the closure of
an annulus $\mathbb A$ is a spectral set for $T$, then $T$ has a
normal boundary of $\mathbb A$-dilation \cite{MR87a:47007}.  Despite a
great deal of effort, there has until recently been little progress in
determining whether or not these examples are typical or exceptional
for finitely connected spectral sets.  The problem may be formulated
in terms of contractive representations of the algebra of rational
functions on the spectral set.  Arveson showed \cite{MR52:15035} that
rational dilation holds precisely when such representations are
automatically completely contractive.  Hence the problem bears some
resemblance to the famous problem of Halmos on similarity to a
contraction solved by Pisier \cite{MR97f:47002}, where the question
hinged on deciding whether bounded representations of the disc algebra
are completely bounded.

Agler, Harland and Rafael \cite{AHR} have recently given an example of
a triply connected domain $R$ and a $4\times 4$ matrix with $X$ as a
spectral set such that this matrix does not have a normal
$B$-dilation.  In this paper we give a proof that for any bounded
triply connected domain with rectifiable disjoint boundary components,
there is an operator $T$ with $X$ as a spectral set which does not
have a normal $B$-dilation.  It can be shown that $T$ can be taken to
be a matrix \cite{VIP}, but the argument does not give a bound on the
size.  Our proof, in contrast to that of Agler, Harland and Raphael,
does not use numerical computations.  Nevertheless it borrows freely
from and overlaps considerably with \cite{AH}, \cite{MR93m:47013},
\cite{AHR}; see also \cite{MR2000a:47034}, \cite{MR2000b:47028},
\cite{MR2002f:47028}.  We implicitly employ the rank one bundle shifts
over $R$ \cite{MR53:1327}, \cite{MR32:6256}, the representation for
the corresponding reproducing kernels in terms of theta functions due
to Fay \cite{MR49:569} and Ball and Clancey \cite{MR97f:46042} (see
also \cite{MR91f:47038}, \cite{MR92e:46110}, \cite{MR2003f:47016},
\cite{MR2000d:47031}), and some elementary compact Riemann surface
theory.

For the remainder of the paper we assume $-1< c_1<0$, $0< c_2<1$ are
two points on the real axis, $0 < r_k < \min{|c_k|,|1-c_k|}$, and that
$R$ is a region obtained by removing disks with centers $c_k$ and
radii $r_k$ from the unit disk.  By \cite{MR37:4243a}, any bounded
triply connected region with boundary components consisting of
disjoint Jordan loops is conformally equivalent to an annulus with a
disk removed.  By scaling and rotating if necessary, we assume that
the outer boundary of the annulus is the unit circle and that the
center of the removed disk is on the real axis.  By choosing a point
between the removed disk and the central disk of the annulus and
applying a M\"obius transformation mapping the unit disk to itself and
the chosen point to the origin, we get a region $R$ of the sort
described initially.  Hence there is no loss of generality in
restricting ourselves to such regions.  Note that in this case, the
boundary $B$ of $R$ consists of three components, $B_0 = \{|z| = 1\}$,
$B_1 = \{|z-c_1| = r_1\}$ and $B_2 = \{|z+c_2| = r_2\}$.  We set $X =
R\cup B$.

As it happens, only three parameters are really needed to distinguish
conformally distinct triply connected domains, so there is some
redundancy in using the four parameters $c_1,c_2,r_1,r_2$.  However,
we will later require that no minimal inner function on our domain
have a zero of multiplicity three at $0$ (all such functions are
normalized to have $0$ as a zero).  We will show that in choosing the
point to move to the origin in going from the annulus with the disk
removed to unit disk with disks on either side of the imaginary axis
removed in the previous paragraph, there is a choice that enforces
this condition on the zeros of these inner functions (Corollary
\ref{cor:notriplezero}).

The main theorem of the paper is the following (see also Agler,
Harland and Rafael \cite{AHR}).

\begin{theorem}
  \label{thm:main}
  For $R$ a bounded triply connected subset of $\mathbb C$ the
  boundary $B$ of which consists of disjoint Jordan curves, there
  exists a Hilbert space $\mathcal H$ and a bounded operator $T$ on
  $\mathcal H$ so that $\overline{R}$ is a spectral set for $T$, but
  $T$ does not have a normal $B$-dilation.
\end{theorem}

Here is the idea of the proof.  We assume without loss of generality
that $R$ is the unit disk with two smaller disks removed, each
centered on the real axis and on opposite sides of the imaginary axis.
Let $\mathcal C$ denote the cone generated by
\begin{equation*}
  \{H(z)(1-\psi(z)\psi(w)^*)H(w)^*:\psi \in B\mathbb H(X), H\in
  M_2(\mathbb H(X))\},
\end{equation*}
where $B\mathbb H(X)$ is the unit ball in the supremum norm of the
space of functions analytic in a neighborhood of $X$, $M_2(\mathbb
H(X))$ the $2\times 2$ matrices of functions analytic in a
neighborhood of $X$.  For $F\in M_2(\mathbb H(X))$, we set
\begin{equation*}
  \rho_F = \sup \{\rho>0:I-\rho^2 F(z)F(w)^* \in \mathcal C\}.
\end{equation*}
It happens that $\rho_F > 0$ and that if $F$ is analytic in a
neighborhood of $X$ and unitary-valued on $B$ such that $\rho_F <1$,
then there exists a Hilbert space $\mathcal H$ and an operator $T$
such that $T$ has $X$ as a spectral set, but $T$ does not have a
normal $B$-dilation (Theorem \ref{thm:possstatz}).  We also show that
if $F\in M_2(\mathbb H(X))$ is as above and $\rho_F =1,$ then under
suitable assumptions regarding zeros, $F$ has a Herglotz
representation (Theorem \ref{thm:rep}).  The kernel of the adjoint of
multiplication by $F$ in this case has its kernel spanned by a finite
collection of reproducing kernels due to Fay (a variant on the Szeg{\H
  o} kernel for harmonic measure---Lemma \ref{lem:rank}), and this
results in $F$ being diagonalizable (Theorem \ref{thm:diagonal}).  Our
counterexample is obtained by constructing a function satisfying the
boundary conditions and hypotheses concerning zeros mentioned above
which is not diagonalizable.  For such a function, $\rho_F <1$, and
hence rational dilation fails.

The organization of the paper differs somewhat from the outline of the
proof mentioned above.  Minimal inner functions play a key role
throughout, so we begin by detailing the results on harmonic functions
and analytic functions with positive real parts which are related to
these inner functions via a Cayley transform.  Much of this material
can be found in Grunsky's monograph \cite{MR57:3365} for more general
planar domains.  These inner functions are used to give a ``scalar''
Herglotz representation theorem for functions analytic in a
neighborhood of $R$ with unimodular boundary values (Proposition
\ref{prop:scalarrep}).  The reader is also referred to \cite{AHR},
where the importance of Herglotz representations in characterizing
spectral sets is cogently presented.

We then seemingly digress into some basic results on Riemann surfaces
and theta functions, which are useful in constructing meromorphic
functions on a compact Riemann surface with a given zero/pole
structure.  The compact surface $Y$ we consider is very special: it is
the double of our two holed region $R$, obtained topologically by
gluing a second copy of $R$ to itself along $B$.  Our minimal inner
functions extend to meromorphic functions on the double via
reflection, and indeed the same goes for matrix valued inner
functions.

Using fairly elementary tools, we are able to say quite a bit about
the zero structure of the minimal inner functions on $R$ based on the
parametrization of these functions.  We then extend results on scalar
inner functions to some $2\times 2$ matrix valued inner functions,
constructing a family of such functions which will ultimately give
rise to our example of such a function with a particular zero
structure which is not diagonalizable.

We then turn to considering Szeg{\H o} kernels $K^a(\zeta,z)$ on $R$
with respect to harmonic measure for the point $a\in R$.  The truly
remarkable fact discovered by Fay \cite{MR49:569} is that these
kernels extend to meromorphic functions on the double Y.  Indeed, Fay
gives a representation for $K^a(\zeta,z)$ in terms of theta functions
and the prime form.  Ball and Clancey \cite{MR97f:46042} give a
similar representation, which is actually a bit more explicit as it
only involves theta functions, for the $n$-torus family of Abrahamse
kernels associated to a multiply connected domain of connectivity
$n+1$.  The explicit nature of the theta function representation of
the meromorphic kernels $K^a(\zeta,z)$ allows the determination of
their zero/pole structure.  Next we consider certain finite linear
combinations of $K^a(\zeta,z)$ with coefficients in ${\mathbb{C}}^2$,
and show that if an analytic function is in the span of certain
kernels, then it must be constant.  This is used later to prove the
diagonalisation result mentioned above.

We next turn to representing nice $2\times 2$ matrix valued inner
functions---these are functions which are unitary valued on the
boundary with what we term a ``standard'' zero set---see section
\ref{subsec:matrixinner} for the definition.  First we show the
connection between the failure of rational dilation and $\rho_F$ being
strictly less than one for some contractive analytic function $F$ with
unitary boundary values.  We make a brief foray into matrix measures,
and then prove that when $\rho_F = 1$ for such functions, we have a
nice representation for $1-F(z)F(w)^*$.  Part of this relies on a
Agler-Nevanlinna-Pick interpolation result (Proposition
\ref{prop:NPinterpolate}), proved using a transfer function
representation, as well as a uniqueness result (Proposition
\ref{prop:unique}).  First we show that a Herglotz type representation
holds over finite subsets of points in $R$ for such an $F$.  This is
then extended to all of $R$ via the interpolation theorem, and if the
finite set of points was chosen in just the right way, this extension
is unique.  The result is what we term a ``tight'' representation.

In the next section we return to Fay's kernel, and use it to prove our
diagonalization result.  We show that some of the matrix inner
functions constructed earlier are not diagonalizable, proving our main
result.

The last section shows why, once we know that rational dilation fails,
that there is a finite dimensional example.  This part is based on
work of Paulsen \cite{VIP}.


\section{Spectral Sets and Some Function Theory on $R$}

We begin with some standard material on spectral sets.  Then we review
harmonic functions, analytic functions of positive real part, and
inner functions on $R$, the unit disk with two disks removed as
described above.

\subsection{Spectral sets}
Let $\mathbb H(X)$ denote complex valued functions which are analytic
in a neighborhood of $X$, and $M_2(\mathbb H(X))$ the $2\times 2$
matrices of such functions.  We likewise let $\mathcal R(X)$ denote
the rational functions with poles off of $X$, with $M_2(\mathcal
R(X))$ defined in the obvious manner.  For $f$ in either $\mathbb
H(X)$ or $M_2(\mathbb H(X))$ define,
\begin{equation*}
 \|f\|=\|f\|_R=\sup\{\|f(z)\|:z\in R\},
\end{equation*}
where $\|f(z)\|$ is the modulus of the scalar $f(z)$ or the operator
norm of the $2\times 2$ matrix $f(z)$ respectively.

Throughout the rest of this section $T$ is a bounded operator on the
complex Hilbert space $\mathcal H$ which has $X$ as a spectral set.

To the operator $T$ associate the homomorphism, $\phi_T:\mathcal
R(X)\to \mathcal B(\mathcal H)$ by $\phi_T(p/q) = p(T)q(T)^{-1}$, $p$
and $q$ polynomials.  When $T$ has a normal $B$-dilation, equation
(\ref{eq:dilation}) can then be expressed as $\phi_T(f)=P_{\mathcal
  H}f(N)|\mathcal H$ for $f\in \mathcal R(X)$.  Using the Riesz
functional calculus and Runge's theorem, the map $\phi_T$ extends
continuously to $\phi_T:\mathbb H(X)\to \mathcal B(\mathcal H)$.
Conversely, a contractive unital homomorphism $\pi:\mathbb H(X)\to
\mathcal B(\mathcal H)$, that is a unital homomorphism satisfying
\begin{equation*}
   \|\pi(f)\|\le \|f\|_R
 \end{equation*}
 for $f\in \mathbb H(X),$ determines an operator with $X$ as a
 spectral set.
 
 If $T$ has a normal $B$-dilation as in (\ref{eq:dilation}) and if $G$
 is in $M_2(\mathcal R(X))$, then
\begin{equation*}
  G(T) = P_{\mathcal H\oplus \mathcal H} G(N) | \mathcal H\oplus
    \mathcal H .
\end{equation*}
Since, by the maximum principle, $\|G(N)\| = \|G\|_R$, it follows that
\begin{equation}
 \label{eq:ssinequality}
  \|G(T)\|\le \|G\|_R.
\end{equation}
Indeed, the same reasoning implies this for any $G\in M_n(\mathcal
R(X))$, or equivalently, that $\phi_T$ is completely contractive.

If $F\in M_2(\mathbb H(X))$ (so that the entries $F_{j,\ell}$ of $F$
are analytic in a neighborhood of $X$), but not necessarily rational,
it then still makes sense to consider, on $\mathcal H\oplus \mathcal
H$, the operator
\begin{equation*}
  F(T) = \begin{pmatrix} F_{j,\ell}(T) \end{pmatrix}
\end{equation*}
where the $F_{j,\ell}(T)$ are defined using the Riesz functional
calculus.

\begin{lemma}
  \label{lem:lem1}
  If $T$ has $X$ as a spectral set and a normal $B$-dilation, and if
  $F\in\mathbb{H}(X)$ is unitary valued on $B$, then
  \begin{equation*}
    \|F(T)\| \le 1.
  \end{equation*}
\end{lemma}

\begin{proof}
  Choose a compact set $K$ so that the interior of $K$ contains $X$
  and $K$ is a subset of the domain of analyticity of $F$.
 
  Using Runge's Theorem (entrywise), there exists a sequence $G_n$ of
  rational $2\times 2$ matrix valued functions with poles off $K$
  which converges uniformly to $F$ on $K$.  From standard results
  about the functional calculus, $\{G_n(T)\}$ converges to $F(T)$ in
  operator norm.  Since $\|F\|_R=1$, the sequence $\{\|G_n\|_R\}$
  converges to $1$.  An application of equation
  (\ref{eq:ssinequality}) and a limit argument completes the proof.
\end{proof}

\begin{remark}
  The set $X$ is a complete spectral set for an operator $T\in
  \mathcal B(\mathcal H)$ if $\|F(T)\|\le \|F\|_R$ for every $n\times
  n$ (no bound on $n$) matrix-valued rational function with poles off
  of $R$ (the norm $\|F\|_R$ defined in the expected way).  It is a
  result of Arveson that if $X$ is a complete spectral set for $T$,
  then $T$ has a normal $B$-dilation \cite{MR40:6274},
  \cite{MR52:15035}, see also \cite{MR88h:46111}, \cite{MR1976867}.
  (Arveson's result is actually stated and proved for a commuting
  $n$-tuple of operators with a domain in $\mathbb C^n$.)
\end{remark}

The following may be found in Conway \cite{MR92h:47026}.

\begin{lemma}
 \label{lem:lem2}
 If $\pi:\mathbb H(X)\to \mathcal B(\mathcal H)$ is a contractive
 unital homomorphism, then $X$ is a spectral set for $T=\pi(\zeta)$.
 Here $\zeta(z)=z$.  Moreover, if $f\in \mathbb H(X)$, then
 $\pi(f)=f(T)$.
\end{lemma}

\begin{proof}
  Given $\lambda \notin X$, the function $f_\lambda(z) =
  (z-\lambda)^{-1}\in \mathbb H(X)$.  Thus, as $\pi$ is a unital
  homomorphism,
 \begin{equation*}
    I=\pi(1)\\
     =\pi(f_\lambda \cdot (\zeta -\lambda))
     =\pi(f_\lambda)\pi(\zeta -\lambda)
     =\pi(f_\lambda)(T-\lambda).
 \end{equation*} 
 It follows that the spectrum of $T$ is in $X$.
 
 If $f=p/q\in \mathcal R(X)$, $p,q$ polynomials, then $\pi(f)=f(T)$
 since $\pi$ is a homomorphism and $f(T)$ can be defined as
 $p(T)q(T)^{-1}$.  Since $\|\pi(f)\|\le \|f\|_R$, it follows that $X$
 is a spectral set for $T$.
 
 Finally, if $f\in H(X)$, then by Runge's Theorem, there exists a
 sequence $\{f_n\}$ from $\mathcal R(X)$ which converges uniformly to
 $f$ on a compact set $K$ containing $X$.  By continuity, both
 $\{f_n(T)\}$ converges to $f(T)$ and $\{\pi(f_n)\}$ converges to
 $\pi(f)$ (in operator norm).  Since $f_n(T)=\pi(f_n)$, it follows
 that $f(T)=\pi(f)$.
\end{proof}

For $F\in M_2(\mathbb H(X))$, define
\begin{equation*} 
  \pi(F)=\begin{pmatrix} \pi(F_{j,\ell}) \end{pmatrix}.
\end{equation*}
Thus, $\pi$ is defined entry-wise.

\begin{proposition}
 \label{prop:extendpi}
 Suppose $\pi:\mathbb H(X)\to \mathcal B(\mathcal H)$ is a contractive
 unital homomorphism and let $T=\pi(\zeta)$.  If $T$ has a normal
 $B$-dilation and if $F$ is an analytic $2\times 2$ matrix valued
 function analytic in a neighborhood of $R$ and unitary valued on $B$,
 then $\|\pi(F)\| \le 1$.
\end{proposition}

\begin{proof}
  By Lemma \ref{lem:lem1}, $\|F(T)\|\le 1$.  On the other hand, by
  Proposition \ref{lem:lem2}, $F(T)=\pi(F)$.
\end{proof}

\subsection{Harmonic functions}

Most of the results and discussion in this section come from Fisher's
book~\cite{MR85d:30001}.

For each point $z\in R$ there exists a measure $\omega_z$ on $B$ such
that if $h$ is continuous on $X$ and harmonic in $R$, then
\begin{equation*}
  h(z) = \int_B h(\zeta) \,d\omega_z(\zeta).
\end{equation*}
The measure $\omega_z$ is a Borel probability measure on $B$ and is
known as \textit{harmonic measure} for the point $z$.

For $z\in R$, $d\omega_z$ and Lebesgue measure $ds$ are mutually
absolutely continuous.  Thus there is a Radon-Nikodym derivative,
\begin{equation*}
  \mathbb P(\cdot ,z) = \frac{d\omega_z}{ds}.
\end{equation*}
This is also the Poisson kernel for $R$; that is,
\begin{equation*}
 \mathbb P(\cdot ,z)=-\frac{1}{2\pi} \frac{\partial}{\partial n}
 g(\cdot ,z),
\end{equation*}
where $g(\cdot ,z)$ is the Green's function for the point $z$, $n$ the
outward normal.

If $h$ is a positive harmonic function in $R$, then there exists a
positive measure $\mu$ such that
\begin{equation}
  \label{eq:repharmonic}
  h(z) = \int_B \mathbb P(\zeta ,z) \,d\mu(\zeta).
\end{equation}
Conversely, given a positive measure $\mu$, the formula
(\ref{eq:repharmonic}) defines a positive harmonic function in $R$.
Note that
\begin{equation*}
  h(0) = \int_B \, d\mu = \mu(B).
\end{equation*}

Let $h_j$ denote harmonic measure for $B_j$.  Thus, $h_j$ is the
solution to the Dirichlet problem with boundary values $1$ on $B_j$
and $0$ on $B_\ell$, $\ell \neq j$.  Alternatively,
\begin{equation*}
    h_j(z) = \int_{B_j} \mathbb P(\zeta ,z) \, ds(\zeta),
\end{equation*}
where $ds$ is normalized arc length measure for $B$ (Nehari
\cite{MR51:13206}, section VII.3).

The symmetry in the domain $R$ yields a simple but useful symmetry for
the $h_j$.

\begin{lemma} 
 \label{lem:propertiesofhj}
 For $z\in R$, $h_j(z)=h_j(z^*)$.
\end{lemma}

\begin{proof}
  This is obvious from the definition of $h_j(z)$.
\end{proof}

\subsection{Analytic functions with positive real part}
As in the last section, some of the results can be found in the book
of Fisher \cite{MR85d:30001}.  The material on extreme points of the
set of normalized analytic functions with positive real part is a
special case of that found in Grunsky \cite{MR57:3365}.

If $\mu$ and $h$ are as in the formula (\ref{eq:repharmonic}), then
the periods $P_j(h)$ of the harmonic conjugate of $h$ around $B_j$, $j
= 0, 1, 2$, are given by
\begin{equation*}
  P_j(h) = \int_B Q_j \,d\mu,
\end{equation*}
where $Q_j$ is the normal derivative of $h_j$.  Of course $h$ is the
real part of an analytic function if and only if $P_j(h) = 0$, $j = 0,
1, 2$.

\begin{lemma}
 \label{lem:zerosofQj}
 The functions $Q_j$ have no zeros on $B$.  Moreover, $Q_j>0$ on $B_j$
 and $Q_j<0$ on $B_\ell$ for $\ell \neq j$.
\end{lemma}

Before proving Lemma \ref{lem:zerosofQj} we note the following
consequence.

\begin{lemma}
 \label{lem:allBj}
 If $h$ is a nonzero positive harmonic function on $R$ which is the
 real part of an analytic function and if $h$ is represented in terms
 of a positive measure $\mu$ as in equation (\ref{eq:repharmonic}),
 then $\mu(B_j)>0$ for each $j$.
\end{lemma}

\begin{proof}
  If $\mu(B_1)=0$, then, as $Q_1<0$ on $B_0 \cup B_2$, $P_1(h)<0$.
  Thus, $\mu(B_1)>0$.  Likewise for $\mu(B_2)$.
 
  Let $h_0$ denote harmonic measure for $B_0$.  Then $\sum_0^2 h_j=1$.
  Consequently, $\sum_0^2 Q_j=0.$ On the other hand, from the proof of
  Lemma \ref{lem:zerosofQj}, $Q_0>0$.  Thus, if $\mu(B_0)=0$, then
  \begin{equation*}
    \begin{split}
      \sum_{j=1}^2 \int_B Q_j \, d\mu =& \sum_{j,\ell=1}^2
      \int_{B_\ell} Q_j \, d\mu \\
      =& -\sum_{\ell=1}^2\int_{B_\ell} Q_0 \, d\mu <0,
    \end{split}
  \end{equation*}
  implying that not both $P_j(h)=0$.  So we must have $\mu(B_0)>0$.
\end{proof}

\begin{proof}[Proof of Lemma \ref{lem:zerosofQj}.]
  Let $h_j$ denote harmonic measure for $B_j$.  Thus, $h_j$ is the
  solution to the Dirichlet problem with boundary values $1$ on $B_j$
  and $0$ on $B_\ell$, $\ell \neq j$.  Alternatively,
  \begin{equation*}
    h_j(z) = \int_{B_j} \mathbb P(\zeta ,z) \, ds(\zeta),
  \end{equation*}
  where $ds$ is normalized arc length measure for $B$.
  
  The functions $Q_j$ are related to the $h_j$ by
  \begin{equation*}
    Q_j = \frac{\partial h_j}{\partial n},
  \end{equation*}
  where the derivative is with respect to the outward normal to the
  boundary.  Note that, for $j = 0, 1, 2$, the partial derivatives
  $\frac{\partial h_j}{\partial n}$ are evidently nonnegative on $B_j$
  and nonpositive on $B_\ell$, for $\ell \neq j$, from which it
  follows that $Q_j\ge 0$ on $B_j$ and $Q_j \le 0$ on $B_\ell$ for
  $\ell \ne j$.
  
  Let $R^\prime$ denote the reflection of $R$ about $B_0$ by $z\mapsto
  1/z^*$.  The functions $h_j$ naturally extend to harmonic functions
  on $X\cup R^\prime$ by
   \begin{equation*}
      h_j(1/z^*) = -h_j(z).
   \end{equation*}
   Thus, $Q_j$ extends to $R^\prime$.
 
   If $Q_j$ were to have infinitely many zeros in $B$, then $Q_j$
   would be zero.  By the Cauchy-Riemann equations, $\tilde{h}_j$, a
   harmonic conjugate of $h_j$ exists on $B_0$ for $j=1,2$.  Moreover,
   since $Q_j\le 0$ on $B_0$ and has only finitely many zeros on
   $B_0$, $\tilde{h}_j$ is strictly decreasing (in the positive
   orientation on $B_0$).  In particular, $\tilde{h}_j$ is one-one on
   $B_0.$ Observe, if $Q_j$ is zero at a point $\zeta \in B_0$, then,
   as the derivative of $h_j$ tangential to $B$ at $\zeta$ is also
   $0$, the function $f=h_j+i\tilde{h_j}$ is analytic near $\zeta$ and
   has zero derivative at $\zeta$, It follows that $f$ is at least two
   to one in a sufficiently small neighborhoods of $\zeta$.  But then
   $\tilde{h_j}$ could not be one to one on $B_0$ near $\zeta$, a
   contradiction.  Hence $Q_j<0$ on $B_0$ ($j=1,2$).
  
   Similar arguments show $Q_j>0$ on $B_j$ and $Q_j<0$ on $B_\ell$ for
   $\ell \ne j$.
\end{proof}

Let $\Pi = B_0\times B_1\times B_2$.

\begin{lemma}
  \label{lem:allpositive}
  For each $p\in\Pi$, the kernel of
  \begin{equation*}
    M(p) = \begin{pmatrix} Q_1 (p_0) & Q_1 (p_1) & Q_1 (p_2) \\ 
      Q_2 (p_0) & Q_2 (p_1) & Q_2 (p_2) \end{pmatrix}
  \end{equation*}
  is one-dimensional and spanned by a vector with all entries strictly
  positive.  In particular, there is a continuous function $\tau:\Pi
  \to \mathbb R^3$ such that $\tau(p)$ is entry-wise positive, the sum
  of the entries is one, and $\tau(p)$ is in the kernel of $M(p)$.
  
  Moreover, $\tau$ reflects the symmetry in the domain.  Namely,
  \begin{equation*}
    \tau(p_0,p_1,p_2)=\tau(p_0^*,p_1,p_2)
     =\tau(p_0,p_1^*,p_2)=\tau(p_0,p_1,p_2^*).
  \end{equation*}
\end{lemma}

\begin{proof}
  Computing the cross product of the rows of $M(p)$, we get the vector
  with entries
  \begin{equation*}
    \begin{split}
      Q_1(p_1)Q_2(p_2)-&Q_1(p_2)Q_2(p_1), \\
      Q_1(p_2)Q_2(p_0)-&Q_1(p_0)Q_2(p_2), \\
      Q_1(p_0)Q_2(p_1)-&Q_1(p_1)Q_2(p_0).
    \end{split}
  \end{equation*}
  This vector is in the kernel of $M(p)$ and by considering the signs
  of the $Q_j$ on the boundary components $B_\ell$, one easily checks
  that the sign of the last two entries are positive.  It is also
  clear from the signs of the entries that the rows of $M(p)$ are
  linearly independent.  Hence $M(p)$ is rank two and its kernel is
  one dimensional.  To finish the proof of the first part of the
  lemma, it remains to show that the sign of the first entry is
  positive.  Since $\sum_{j=0}^2 h_j(P)= 1$ for $P\in X$, it follows
  that for $P\in B_1$ or $P\in B_2$, that
  \begin{equation*}
    Q_1(P)+Q_2(P)=-Q_0(P)>0.
  \end{equation*}
  Thus,
   \begin{equation*}
    \begin{split}
      Q_1(p_1)Q_2(p_2)-& Q_1(p_2)Q_2(p_1)\\
      =& Q_1(p_1)Q_2(p_2) + Q_1(p_1)Q_1(p_2) - Q_1(p_1)Q_1(p_2)-
      Q_1(p_2)Q_2(p_1)\\
      =& Q_1(p_1)(Q_2(p_2)+Q_1(p_2))- Q_1(p_2)(Q_1(p_1)+Q_2(p_1))\\
      =& - Q_1(p_1) Q_0(p_2) +Q_1(p_2)Q_0(p_1).
    \end{split}
  \end{equation*}
  Examination of the signs of the terms on the right hand side above
  shows that the first term is indeed positive.
  
  To prove the last part of the lemma, simply note that the symmetry
  in the domain implies each $h_j$ is symmetric, $h_j(z^*)=h_j(z)$ and
  thus, $Q_j(z^*)=Q_j(z)$.
\end{proof}

The lemma allows the construction of canonical analytic functions of
positive real part on $R$, since, by construction,
\begin{equation}
 \label{eq:definehp}
  h_p = \sum \tau_j(p) \, \mathbb P(\cdot,p_j)
\end{equation}
is a positive harmonic function with no periods.  Indeed, $h_p$
corresponds to the measure $\sum \tau_j(p)\, \delta_{p_j}$ in equation
(\ref{eq:repharmonic}).  In particular, by our normalization of
$\tau$,
\begin{equation*}
  h_p(0) = 1.
\end{equation*}
Let $f_p$ denote the analytic function such that $f(p) = 1$ and with
real part $h_p$.

Let $H(R)$ be the locally convex metrizable topological space of
holomorphic functions on $R$ with the topology of uniform convergence
on compact subsets of $R$.  The space $H(R)$ has the Heine-Borel
property; i.e., closed bounded subsets of $H(R)$ are compact.  Let
\begin{equation*}
  \mathbb K = \{f\in H(R): f(0) = 1, f+f^*>0\}.
\end{equation*}
The set $\mathbb K$ is easily seen to be closed.

\begin{lemma} 
 \label{lem:Kcompact}
 The set $\mathbb K$ is compact.
\end{lemma}

What needs to be shown is that $\mathbb K$ is bounded---in other words
that for each compact subset $K\subset R$ there exists an $M_K$ so
that $|f(z)|\le M_K$ for all $f\in\mathbb K$ and $z\in K$.

That $\mathbb K$ is bounded is straightforward in the case that $R$ is
replaced by the unit disk, and we consider this case to begin with.

Let $P_r(\theta)$ denote the Poisson kernel and $Q_r(\theta)$ the
conjugate Poisson kernel,
\begin{equation*}
 \begin{split}
   P_r(\theta) = &\frac{1-r^2}{1+r^2-2r\cos(\theta)} \\
   Q_r(\theta) = &\frac{2r\sin(\theta)}{1+r^2-2r\cos(\theta)}.
 \end{split}
\end{equation*}
If $U$ is a positive harmonic in the unit disk $\mathbb D$ with
$U(0)=1$, then there is a probability measure $\mu$ so that
\begin{equation*}
 U(r\exp(i\theta))=\int_{-\pi}^\pi P_r(\theta-t)\, d\mu(t).
\end{equation*}
The harmonic conjugate $V$ of $U$ normalized by $V(0)=0$ is given by
\begin{equation*} 
 V(r\exp(i\theta))=\int_{-\pi}^\pi Q_r(\theta-t)\, d\mu(t).
\end{equation*}
Thus, as $\mu$ is a probability measure, for each $0<r<1$ there is a
constant $M_r$ so that if $|z|\le r$, then $|U(z)|,|V(z)|\le M_r/2$.
By conformal mapping, the result is seen to hold for any bounded
simply connected domain.

The proof of the lemma is based upon the above result for simply
connected domains, using the fact that $R$ can be written as the union
of two simply connected domains.

\begin{proof}[Proof of Lemma \ref{lem:Kcompact}.]
  It suffices to show that for each compact set $K$ there is a
  constant $M_K$ so that if $f\in \mathbb K$ and $z\in K$ so that
  $|f(z)|\le M_K$.  In fact, it is enough to work with the compact
  sets of the form
  \begin{equation*}
    K_\epsilon =\{|z|\leq 1-\epsilon\}\cap\{|z+c_1|\geq r_1+\epsilon\}
    \cap \{|z+c_2| \geq r_2+\epsilon\},
 \end{equation*}
 for $\epsilon > 0$ sufficiently small.
  
 First, observe that if $f$ is holomorphic with positive real part on
 $R$ and $f(0)=1$, then for each $\epsilon>0$ there is a constant
 $N_\epsilon$ so that for all $z\in K_\epsilon$, $h(z)\le N_\epsilon$,
 by virtue of the representation (\ref{eq:repharmonic}).
  
 Given $\epsilon>0$, the open simply connected sets
  \begin{equation*}
    R^\pm_\epsilon=\left\{ z=a+bi \in K_{\frac{\epsilon}{2}}: \pm b\ge
      - \frac{\min(r_1,r_2)}{2}\right\}
  \end{equation*}
  contain $0$ and their union is $K_{\frac{\epsilon}{2}}$.  After
  conformal mapping from the disk taking $0$ to $0$, it follows that
  there exists $M_{\epsilon}^\pm$ so that if $f$ is analytic on $R$
  with $f(0)=1$, then for all $z\in R^{\pm}_\epsilon$ $|f(z)|\le
  M_{\epsilon}^\pm$.
\end{proof}

\begin{lemma}[\cite{MR57:3365}]
  \label{lem:extremepoints}
  The extreme points of $\mathbb K$ are precisely $\{f_p:p\in \Pi\}.$
\end{lemma}

\begin{proof}
  It is evident that each $f_p$ is an extreme point of $\mathbb K$.
  So consider the converse.
  
  Let $f\in \mathbb K$.  The real part of $f$ is a positive harmonic
  function $h$ with $h(0) = 1$.  Hence there exists a probability
  measure $\mu$ such that
  \begin{equation*}
    h(z) = \int_B \mathbb P(\zeta ,z)\,d\mu(\zeta).
  \end{equation*}
  
  Suppose the support of $\mu$ on $B_0$ contains more than one point
  and so can be written as the union of disjoint sets $A_1, A_2
  \subset B_0$ with $\mu(A_j)>0$.
  
  Let
  \begin{equation*}
    \alpha_{j,\ell} = \int_{A_\ell} Q_j \,d\mu , \qquad \ell = 1,2.
  \end{equation*}
  and
  \begin{equation*}
    \kappa_{j,m} = \int_{B_m} Q_j \,d\mu , \qquad m = 1,2.
  \end{equation*}  
  Since $h$ is the real part of an analytic function
  \begin{equation*}
    0 = \int_B Q_j \,d\mu.
  \end{equation*}
  Thus, $\kappa_{j,1} + \kappa_{j,2} +(\alpha_{j,1} +\alpha_{j,2})
  =0$.  Since also $Q_j<0$ on $B_0$ for $j=1,2$,
  \begin{equation*}
    \kappa_{j,1} + \kappa_{j,2} = -(\alpha_{j,1} +\alpha_{j,2}) > 0.
  \end{equation*}
  This gives, $\kappa_{1,1}\ge |\kappa_{1,2}|=-\kappa_{1,2}$.  Hence
  the the determinant of $\kappa = (\kappa_{i,j})$ is positive.
  
  Since the determinant of $\kappa =(\kappa_{i,j})$ is positive the
  solution of
  \begin{equation*}
    \begin{pmatrix}
      \kappa_{1,1} & \kappa_{1,2} \\ \kappa_{2,1} & \kappa_{2,2}
    \end{pmatrix}
    \begin{pmatrix}
      \beta_{1,\ell} \\ \beta_{2,\ell}
    \end{pmatrix}
    = -
    \begin{pmatrix}
      \alpha_{1,\ell} \\ \alpha_{2,\ell}
    \end{pmatrix}
  \end{equation*}
  is given by
  \begin{equation*}
   \frac{1}{\det{k}} \begin{pmatrix}
      \kappa_{2,2} & -\kappa_{1,2} \\ -\kappa_{2,1} & \kappa_{1,1}
    \end{pmatrix}
    \begin{pmatrix}
      -\alpha_{1,\ell} \\ -\alpha_{2,\ell}
    \end{pmatrix}
    = 
    \begin{pmatrix}
      \beta_{1,\ell} \\ \beta_{2,\ell}
    \end{pmatrix} .
  \end{equation*}
  In view of the signs of the $\kappa_{j,m}$ and $\alpha_{j,\ell}<0$,
  it follows $\beta_{1,\ell}, \beta_{2,\ell} \geq 0$.
  
  Define positive measures $\nu_\ell$ by
  \begin{equation*}
    \nu_\ell(A) = \mu(A\cap A_\ell)+ \beta_{1,\ell}\, \mu(A\cap B_1) +
    \beta_{2,\ell}\, \mu(A\cap B_2).
  \end{equation*}
  Then
  \begin{equation*}
    \int_B Q_j \,d\nu_\ell = \alpha_{j,\ell} +
    \kappa_{j,1}\, \beta_{1,\ell}+  \kappa_{j,2}\, \beta_{2,\ell} = 0
  \end{equation*}
  and therefore each
  \begin{equation*}
    h_\ell = \int_B \mathbb P(\cdot,\zeta)\,d\nu_\ell(\zeta)
  \end{equation*}
  is the real part of an analytic function $g_\ell$ with $g_\ell(0) =
  \nu(B)$.  Since also $\nu_1+\nu_2 = \mu$, we have $h_1+h_2 = h$.
  Thus, $\frac{g_\ell}{g_\ell(0)} \in \mathbb K$ and
  \begin{equation*}
    f = g_1(0) \left( \frac{g_1}{g_1(0)}\right) 
    + g_2(0) \left( \frac{g_2}{g_2(0)}\right).
  \end{equation*}
  We conclude that $f$ is not an extreme point.
  
  Next suppose $B_1$ is a disjoint union of sets $A_1,A_2$ with
  $\mu(A_j)>0$.  Let
    \begin{equation*}
    \alpha_{j,\ell} = \int_{A_\ell} Q_j \,d\mu , \qquad \ell = 1,2,
  \end{equation*}
  and
  \begin{equation*}
    \kappa_{j,m} = \int_{B_m} Q_j \,d\mu , \qquad m = 0,1.
  \end{equation*}  
  This time the signs are $\alpha_{1,\ell}\ge 0$, $\alpha_{2,\ell}\le
  0$, $\kappa_{j,0}\le 0$, $\kappa_{1,2}<0$, and $\kappa_{2,2}>0$.
  So, with
   \begin{equation*}
     \kappa=\begin{pmatrix}
             \kappa_{1,0} & \kappa_{1,2} \\ \kappa_{2,0} & \kappa_{2,2}
             \end{pmatrix}
   \end{equation*}
   the determinant of $\kappa$ is negative.  Thus, the solution of
  \begin{equation*}
   \kappa  \begin{pmatrix}
      \beta_{0,\ell} \\ \beta_{2,\ell}
    \end{pmatrix}
    = -
    \begin{pmatrix}
      \alpha_{1,\ell} \\ \alpha_{2,\ell}
    \end{pmatrix}
  \end{equation*}
  is given by
  \begin{equation*}
   \frac{1}{\det{k}} \begin{pmatrix}
      \kappa_{2,2} & -\kappa_{1,2} \\ -\kappa_{2,0} & \kappa_{1,0}
    \end{pmatrix}
    \begin{pmatrix}
      -\alpha_{1,\ell} \\ -\alpha_{2,\ell}
    \end{pmatrix}
    = 
    \begin{pmatrix}
      \beta_{1,\ell} \\ \beta_{2,\ell}.
    \end{pmatrix}
  \end{equation*}
  In view of the signs of the entries and of $\det(\kappa)$, the
  $\beta_{j,\ell}$ are all nonnegative.
  
  Define positive measures $\nu_\ell$ by
  \begin{equation*}
    \nu_\ell(A) = \mu(A\cap A_\ell)+ \beta_{0,\ell}\, \mu(A\cap B_0) +
    \beta_{2,\ell}\, \mu(A\cap B_2).
  \end{equation*}
  Then
  \begin{equation*}
    \int_B Q_j \,d\nu_\ell = \alpha_{j,\ell} +
    \kappa_{j,0} \, \beta_{0,\ell} +  \kappa_{j,2}\, \beta_{2,\ell} = 0
  \end{equation*}
  and therefore each
  \begin{equation*}
    h_\ell = \int_B \mathbb P(\cdot,\zeta)\,d\nu_\ell(\zeta)
  \end{equation*}
  is the real part of an analytic function $g_\ell$ with $g_\ell(0) =
  \nu(B)$.  Since also $\nu_1+\nu_2 = \mu$, we have $h_1+h_2 = h$.
  The argument proceeds as before, with $\frac{g_\ell}{g_\ell(0)} \in
  \mathbb K$ and
  \begin{equation*}
    f = g_1(0) \left( \frac{g_1}{g_1(0)}\right) 
    + g_2(0) \left( \frac{g_2}{g_2(0)}\right).
  \end{equation*}
  We conclude $f$ is not an extreme point.
\end{proof}

\begin{lemma}
  \label{lem:extremeclosed}
  The set of extreme points of $\mathbb K$ is a closed set and the
  function taking $\Pi$ to $\mathbb K$ by $p \mapsto f_p$ is a
  homeomorphism.
\end{lemma}

\begin{proof}
  It suffices to show, if $p(n) \in \Pi$ converges to $p(0)$ in $\Pi$,
  then $f_{p(n)}$ converges uniformly on compact subsets of $R$ to
  $f_{p(0)}$.  Since $\mathbb K$ is compact and $\{f_{p(n)}\}$ is a
  sequence in $\mathbb K$, some subsequence, still denoted
  $\{f_{p(n)}\}$, converges to some $f \in \mathbb K$.  Let $h_{p(n)}$
  denote the real part of $f_{p(n)}$ and let $\mu_n$ denote the
  measure which represents $h_{p(n)}$ so that
  \begin{equation*}
    h_{p(n)}(z) = \int_B \mathbb P(\zeta ,z)\,d\mu_n
  \end{equation*}
  and $\mu_n=\sum \tau_j (p(n))\, \delta_{p(n)_j}$, where $\tau$ is
  defined in Lemma \ref{lem:allpositive}.  See also equation
  (\ref{eq:definehp}).
  
  The measures $\mu_n$ converge to the measure $\mu_0$ weakly so that
  $h_{p(n)}(z)$ converges to $h_{p(0)}(z)$ pointwise in $R$.  It
  follows that $h_{p(0)}$ is the real part of $f$.  Thus,
  $\{f_{p(n)}\}$ converges, in $H(R)$, to $f_{p(0)}$.  Since every
  subsequence of our original subsequence has a subsequence which
  converges to $f_{p(0)}$, the whole subsequence converges to
  $f_{p(0)}$.  This shows that the mapping is continuous.  Since $\Pi$
  is compact, it follows that our mapping is a homeomorphism and its
  range is compact.
\end{proof}

\subsection{Scalar inner functions on $R$}

Up to post composition by a M\"obius transformation, the inner
functions on $R$ with precisely three zeros are canonically
parameterized by the $2$-torus, $\mathbb T^2$.  Details may be found
in the book by Fay \cite{MR49:569}.  This section contains an
alternate description of this family using results from the previous
sections.  The following is well known (see, for example, Fisher
\cite{MR85d:30001}, Ch.~4, ex.~6,7).

\begin{proposition}
  \label{prop:three-zeros}
  A nonconstant inner function $\psi$ on $R$ has at least three zeros
  counting with multiplicity.  Moreover, if $\psi$ has exactly three
  zeros, $z_0,z_1,z_2$, then for $j=1,2$
  \begin{equation*}
    \sum_{\ell=0}^2  h_j(z_\ell)=1.
  \end{equation*}
\end{proposition}

As a first application, the last proposition allows us to show that we
may assume that our region supports no inner function with a zero of
multiplicity three at $0$, and so consequently we will take it that
$R$ has this property throughout the remainder of the paper.

\begin{corollary}
  \label{cor:notriplezero}
  We may assume without loss of generality that there is no inner
  function on $R$ with a zero of multiplicity three at $0$.
\end{corollary}

\begin{proof}
  Suppose $f$ is an inner function on $R$ with exactly three zeros
  counting multiplicity.  Then if $f$ takes the value $f(x)$ with
  multiplicity three where $x\in (c_1+r_1, c_2-r_2)$ (ie, $f(z)-f(x)$
  has a zero of multiplicity three at $x$), then there is a M\"obius
  transformation $m$ of the disk to itself fixing $1$ and moving
  $f(x)$ to $0$.  Thus $m\circ f$ is an inner function with zero of
  multiplicity three at the origin.  Hence by Proposition
  \ref{prop:three-zeros}, $h_j(x) = 1/3$.  Since $h_1$ is continuous
  and equal to $1$ on $B_1$, if we choose $x$ close enough to
  $c_1+r_1\in B_1$, it will be the case that $h_1(x) > 1/3$, and so
  $f$ does not take the value $f(x)$ with multiplicity three for any
  inner function with three zeros.  Now take a M\"obius transformation
  mapping $x$ to the origin.  This takes $R$ to a new region which
  again has the property that it is a region which is the unit disk
  with two smaller disks removed on either side of the imaginary axis
  and centers on the real axis.
\end{proof}

For $p\in \Pi$, let
\begin{equation*}
  \phi_p = \frac{f_p-1}{f_p+1},
\end{equation*}
where $f_p$ is the extreme point of $\mathbb K$ corresponding to the
point $p = (p_0,p_1,p_2)$ in $\Pi$ as in the previous subsection.  The
real part, $h_p$, of $f_p$ is harmonic across $B\setminus
\{p_0,p_1,p_2\}$ and therefore $f_p$ is (at least locally) analytic
across $B\setminus \{p_0,p_1,p_2\}$.  Further, near $p_j$, $f_p$ has
the form $\frac{g_j(z)}{z-p_j}$ for some $g_j$ analytic in a
neighborhood of $s_j$ and nonvanishing at $p_j$ (\cite{MR85d:30001},
Ch.~4).  From these properties of $f_p$, it follows that $\phi_p$ is
continuous onto $B$ and $|\phi_p| = 1$ on $B$.  By the reflection
principle, $\phi_p$ extends to be analytic in a neighborhood of $B$.
Since $\phi_p$ is inner and extends analytically across $B$, and
$\phi_p^{-1}(\{1\}) = \{p_0,p_1,p_2\}$, it follows that preimage of
each point $z\in\mathbb D$ is exactly three points, counted with
multiplicity.  In particular, $\phi_p$ has precisely three zeros.

On the other hand, suppose $\phi$ is analytic in a neighborhood of
$R$, has modulus one on $B$, and three zeros in $R$.  As above, it
follows that $\phi^{-1}(1)$ consists of three points.  Moreover, the
real part of
\begin{equation*}
  f = \frac{1+\phi}{1-\phi}
\end{equation*}
is a positive harmonic function which is zero on $B$ except at those
points $z$ where $\phi(z) = 1$.  By Lemma \ref{lem:allBj} we must have
$\phi^{-1}(\{1\}) = \{p_0,p_1,p_2\}$ where $p_j\in B_j$.  So if we
also assume $\phi(0) = 0$, then $\phi = \phi_p$ for some $p\in\Pi$.

\begin{proposition}
  \label{prop:scalarrep}
  If $\psi$ is analytic in analytic in $R$ and if $|\psi|\le 1$ on
  $R$, then there exists a positive measure $\mu$ on $\Pi$ and a
  measurable function $h$ defined on $\Pi$ whose values are functions
  $h(\cdot,p)$ analytic in $R$ so that
  \begin{equation*}
    1-\psi(z)\psi(w)^* = \int h(z,p)[1-\phi_p(z)\phi_p(w)^*] h(w,p)^*
    \,d\mu(p).
  \end{equation*}
\end{proposition}

\begin{proof}
  First suppose $\psi(0) = 0$.
  
  Let
  \begin{equation*}
    f = \frac{1+\psi}{1-\psi}.
  \end{equation*}
  Verify
   \begin{equation*}
    \psi =\frac{f-1}{f+1}
   \end{equation*}
   and hence,
  \begin{equation}
   \label{eq:scalarrep1}
     1-\psi(z)\psi(w)^* = 2 \frac{f(z)+f(w)^*}{(f(z)+1)(f(w)^*+1)}.
   \end{equation}
   
   Since $h$, the real part of $f$, is positive and $f(0) = 1$, the
   function $f$ is in $\mathbb K$ defined in the previous subsection.
  
   Since $\mathbb K$ is a compact subset of the topological vector
   space $H(R)$ and the extreme points $\{f_p:p \in \Pi\}$ of $\mathbb
   K$ is a compact set by Lemma \ref{lem:extremeclosed}, there exists
   a (regular Borel) probability measure $\nu$ on $\Pi$ so that
  \begin{equation*}
    f = \int_{\Pi} f_p \,d\nu(p).
  \end{equation*}
  
  Using the definition of $\phi_p$ and equation (\ref{eq:scalarrep1}),
  verify,
  \begin{equation*}
    1-\psi(z)\psi(w)^* = \int_{\Pi} \frac{1-\phi_p(z)\phi_p(w)^* }
     {(f(z)+1)(1-\phi_p(z))(1-\phi_p(w)^*)(f(w)^*+1)}\,d\nu(p). 
  \end{equation*}
  
  If $\psi(0) = a$, then one has a representation as above since
  \begin{equation*}
    1-\left( \frac{\psi(z)-a}{1-a^*\psi(z)}\right) 
    {\left( \frac{\psi(w)-a}{1-a^*\psi(w)}\right)}^*  =
    \frac{(1-aa^*)(1-\psi(z)\psi(w)^*)}{(1-a^*\psi(z))(1-a\psi(w)^*)}.
  \end{equation*}
\end{proof}

While we have used the three parameters in $\Pi$ to parameterize the
inner functions with exactly three zeros, after rotation, really only
two are needed.  Indeed,
\begin{equation*}
  \phi_p(1)^* \phi_p = \phi_q,
\end{equation*}
where $q = (1,q_1,q_2)$ and $q_j\in B_j$ for $j=1,2$, are the unique
points such that $\phi_p(q_j)=\phi_p(1)$.


\section{Some Riemann Surfaces and Theta Functions}
We review some of the results from the theory of Riemann surfaces
which we will need subsequently.  In particular, we look at theta
functions, the use of which in operator theory was pioneered by
Clancey \cite{MR91f:47038}, \cite{MR92e:46110}.  The presentation here
borrows heavily from Ball and Clancey \cite{MR97f:46042} as well as
Mumford \cite{MR85h:14026}, \cite{MR86b:14017}, and Farkas and Kra
\cite{MR93a:30047}.

Let $Y$ denote the double of the bordered Riemann surface $X=R\cup B$.
Recall $Y$ is obtained topologically by gluing a second copy
$R^\prime$ of $R$ along $B$.  The complex atlas is then found by
anti-holomorphically reflecting the complex structure from $X$ to
$R^\prime$.  Since $B_0$ is the unit circle, the anti-copy $R^\prime$
of $R$ may be thought of as the reflection of $R$ in $B_0$ given by
$z\mapsto z^{*-1}$.  In particular, there is an anti-holomorphic
involution $J:Y\to Y$ which fixes $B$.  For $\zeta \in R$, $J\zeta$ is
its twin in $R^\prime$.

Note that $Y$ is a compact Riemann surface.  We are interested in it
since we can sometimes extend analytic functions defined on $R$ to
meromorphic functions on $Y$.  In particular this is true for inner
functions, which extend by reflection.

We get a lot of mileage out of a few basic tools for analyzing
meromorphic functions on compact Riemann surfaces.  To begin with, if
such a function is nonconstant and takes a given value $n$ times, it
takes all of its values $n$ times (counting multiplicity).
Furthermore, and as a consequence, such a function has equal numbers
of zeros and poles.  And finally, a meromorphic function on a compact
Riemann surface without zeros or poles is constant.

The material in subsections
\ref{subsec:period-matrix-abel}--\ref{subsec:prime-form} is solely for
describing the zeros of the Fay kernel in
Section~\ref{sec:fay-kernels}.  The reader who is willing to accept
the statement of Theorem \ref{thm:FaysK} may skip these parts.

\subsection{Minimal Meromorphic Functions on $R$}
\label{subsec:minimal}
Our Riemann surface $Y$ also comes naturally equipped with the a
conformal involution which fixes $6$ points, namely
\begin{equation*}
  \iota(\zeta)=J(\zeta^*).
\end{equation*}
Geometrically, $\iota$ is rotation by $\pi$ in the axis through $\pm
1$.  The fixed points of the involution, are precisely the Weierstrass
points of $Y$, namely $\pm 1, c_1\pm r_1, c_2\pm r_2$
(\cite{MR93a:30047}, Cor.~1, p.~108).

If $f$ is a meromorphic function on a compact Riemann surface, then
there is a number $n$, the degree of $f$, so that $f$ takes each value
$n$ times, counting multiplicity.  In particular, if $f$ has just one
pole, then $f$ is one to one and the Riemann surface is conformally
equivalent to the Riemann sphere.  Thus, a nonconstant meromorphic
function on $Y$ must have at least two poles (and zeros).

Up to post composition by a M\"obius transformation, there is a unique
meromorphic function $\phi$ on $Y$ with precisely two poles.
Moreover, $\phi$ is ramified with branching one at the Weierstrass
points.

The construction of a meromorphic function with two poles and zeros
found in Farkas and Kra (Theorem III.7.3) is the following.  The
Riemann surface $\tilde{Y}$, obtained as the quotient of $Y$ by the
map $\zeta\mapsto \iota^{-1}(\zeta)$, is a Riemann surface of genus
zero by the Riemann-Hurwitz formula.  Thus, $\tilde{Y}$ is conformally
equivalent to the Riemann sphere and the canonical quotient map
\begin{equation*}
 Y\mapsto \tilde{Y} 
\end{equation*}
followed by a conformal map from $\tilde{Y}$ to the Riemann sphere
generates the desired meromorphic function.  From this construction it
is evident that if $\phi$ is any meromorphic function with two poles
and two zeros, then $\phi \iota =\phi$.  Another way to see this is
that $\phi - \phi\iota$ has zeros at the six Weierstrass points, while
it has at most four poles, and so must be constantly zero.

\begin{proposition}
 \label{prop:minimalpoles}
 If $f$ is a meromorphic function on $Y$ with precisely two poles
 $Q_1$ and $Q_2$, then $Q_2=\iota(Q_1)$.  In particular, $f$ cannot
 have a pole in $R$ and the other in $X = R\cup B$.  The same
 statement holds if we consider zeros rather than poles.
\end{proposition}

\begin{proof}
  From the discussion above, $f\circ \iota=f$.  The conclusion for the
  poles of $f$ follows.  By considering $1/f$ the conclusion for zeros
  of $f$ is also seen to hold.
\end{proof}

Using the proposition it is possible to show that a unimodular
function on $R$ which takes the value $1$ at exactly one point on each
boundary component $B_j$ and is $0$ at $0$ is uniquely determined.  Of
course, this is also evident by construction and has already been
used.

\begin{corollary}
  \label{cor:uniqueinner}
  Suppose $p_j\in B_j$, for $j=0,1,2$.  If $f,g:X\to \mathbb C$ are
  analytic in a neighborhood of $X$, are unimodular on $B$,
  $f^{-1}(\{1\})\cap X=g^{-1}(\{1\})\cap X=\{p_0,p_1,p_2\}$, and
  $f(0)=g(0)=0$.  Then $f=g$.
\end{corollary}

\begin{proof}
  By the maximum modulus principle, $f$ and $g$ are strictly less than
  one in modulus in $R$.  The functions $f$ and $g$ reflect to
  meromorphic functions on $Y$, which can take the value $1$ only on
  the boundary of $B$.  Thus, the functions are globally three to one
  and so have three zeros in $R$ and three poles, at the reflected
  points, in $R^\prime$.  In particular, $f$ and $g$ share the zero
  $0$ and pole $J0$.
  
  Consider the function
 \begin{equation*}
   \psi=\frac{1-f}{1-g}.
 \end{equation*}
 The zeros of $1-f$ and those of $1-g$ agree, so that $\psi$ can only
 have poles at the poles of $f$.  Since one of these poles matches a
 pole of $g$, $\psi$ has only two poles (and two zeros).  However,
 these poles are poles of $f$ which both are in $R^\prime$,
 contradicting the proposition.  Hence, it must be that $\psi$ is
 constant so that $1-f=c(1-g)$.  Since $f(0)=g(0)=0$, $c=1$.
\end{proof}

We can also describe minimal inner functions which agree at two of
their zeros.

\begin{lemma}
  \label{lem:twozeros}
  If $\phi_p$ and $\phi_{\tilde p}$ have two zeros in common ($0$ and
  one other) and $\phi_p(1) = \tilde\phi_p (1) = 1$, then either they
  are equal or the remaining zeros form a conjugate pair.  In case
  that $p=(1,p_1,p_2)$ and $\tilde p = (1,p_1,{\tilde p}_2)$, then
  ${\tilde p}_2 = p_2$ and so $\phi_{\tilde p} = \phi_p$.
\end{lemma}

\begin{proof}
  If $\phi_p$ and $\phi_{\tilde p}$ have all three zeros in common,
  then they are equal, and there is nothing to prove.  So assume
  instead that they have two common zeros ($0$ and $z_1$) and unequal
  zeros $z_2$, ${\tilde z}_2$ for $\phi_p$ and $\phi_{\tilde p}$,
  respectively.  Then $\phi_p/\phi_{\tilde p}$ has two simple zeros
  and poles, $z_2$, $J{\tilde z}_2$ and ${\tilde z}_2$, $Jz_2$,
  respectively.  By Proposition \ref{prop:minimalpoles}, $Jz_2 = \iota
  {\tilde z}_2 = J({\tilde z}^*_2)$, or ${\tilde z}_2 = z_2^*$.
  
  Now suppose that $\phi_p$ and $\phi_{\tilde p}$ have common zeros
  $0$ and $z_1$, and that $p=(1,p_1,p_2)$ and $\tilde p =
  (1,p_1,{\tilde p}_2)$.  Then unless $\phi_p/\phi_{\tilde p}$ is
  constant, $\phi_p/\phi_{\tilde p} - 1$ has the same two poles as
  $\phi_p/\phi_{\tilde p}$.  It would then have exactly two zeros,
  namely $1$ and $p_1$.  Now apply Proposition \ref{prop:minimalpoles}
  and note that $\iota(1) = 1\neq p_1$.  The result then follows.
\end{proof}

\begin{lemma}
  \label{lem:real-zeros}
  Let $p=(1,p_1,p_2)$.  If the zeros of $\phi_p$ are all real, then
  $p_1,p_2\in\mathbb{R}$.
  
  Conversely, for $p=(1,c_1-r_1,c_2-r_2)$, the zeros of $\phi_p$ are
  are all real.  In fact, $\phi_p$ has a zero at zero and one zero in
  each of the intervals $(-1,c_1-r_1)$ and $(c_2+r_2,1)$ on the real
  axis in $R$.
\end{lemma}

It turns out that for $p=(1,c_1+r_1,c_2-r_2)$ or
$p=(1,c_1-r_1,c_2+r_2)$, the zeros of $\phi_p$ are also real.

\begin{proof}
  Suppose that the zeros of $\phi_p$ are real, and define
  ${\tilde\phi}_p = \phi_p(\zeta^*)^*$.  Then ${\tilde\phi}_p$ has the
  same zeros as $\phi_p$.  Hence the ratio of these two functions has
  no zeros and so must be constant.  Since both equal $1$ at $1$, it
  follows that they are equal everywhere.  In particular,
  \begin{equation*}
    \phi_p(c_j\pm r_j)^* = \phi_p((c_j \pm r_j)^*)^*
    = \phi_p(c_j\pm r_j)
  \end{equation*}
  so that $\phi_p(c_j\pm r_j)$ is real (in each of the four cases).
  Since $\phi_p$ maps each $B_j$ one-one onto the unit circle, it
  follows that that $\phi_p$ takes the value $1$ at one of the points
  $c_j\pm r_j$ on $B_j$, $j=1,2$ and this gives the conclusion of the
  lemma.
      
  For the converse, define $\tilde{\phi}_p$ as above.  Since the
  entries of $p$ are real, Corollary \ref{cor:uniqueinner} implies
  $\phi_p(\zeta)=\tilde{\phi}_p(\zeta^*)^*$ so that $\phi_p$ is real
  on the real axis.  Since also $\phi_p$ takes the value $1$ at
  $c_j-r_j$ and by assumption there is only one point on $B_j$ where
  $\phi_p$ equals $1$, it takes the value $-1$ at the points
  $c_j+r_j$.  The intermediate value theorem now implies that $\phi_p$
  has at least one zero in each of the intervals, $(-1,c_1-r_1)$,
  $(c_1+r_1,c_2-r_2)$, and $(c_2+r_2,1)$.  Since $\phi_p$ has exactly
  three zeros the result follows.
\end{proof}

For $t = (t_1,t_2)$, $t_1,t_2\in [0,2\pi)$, let $p_1(t)=
c_1-r_1\exp(it_1)$ and $p_2 = c_2-r_2\exp(it_2)$.  These points are on
the two inner components of the boundary.  On the outer boundary we
just use $1$.  Accordingly, let $p(t) = (1,p_1(t),p_2(t)) \in\Pi$ and,
for notational ease, let $\psi_t=\phi_{p(t)}$.  Then $\psi_{-t} =
\phi_{p(-t)}$.

\begin{lemma}
 \label{lem:psit}
 There exists $t_1, t_2$ nonzero such that $\psi_t$ has three distinct
 zeros $0, z_1, z_2$ with $z_1, z_2$ both nonreal.
\end{lemma}

\begin{proof}
  Write $0$ for $(0,0)$.  From the second part of Lemma
  \ref{lem:real-zeros}, $\psi_0$ has distinct real zeros.
  
  Now suppose that $t=(t_1,0)$.  If $\pi>t_1>0$ is small enough, the
  zeros $0,z_1,z_2$ of $\psi_t$ are still distinct.  By the first part
  of Lemma \ref{lem:real-zeros}, if these zeros are all real, then
  $t_1$ is a multiple of $\pi$.  On the other hand suppose that one of
  the nonzero zeros of $\psi_t$, say $z_1$, is real.  Then $\psi_t$
  and $\psi_{-t}$ share $0$ and $z_1$ as zeros, as well as both being
  equal to $1$ at $1$ and $c_2-r_2$.  Hence by Lemma
  \ref{lem:twozeros} they are equal, implying that $p_1(t)$ is real,
  or equivalently, that $t_1$ is a multiple of $\pi$.  Since we have
  chosen $0<t_1<\pi$, this cannot happen.  Thus we must have both
  zeros of $\psi_t$ are nonreal for small enough $t_1$.

  Finally, if we choose $t_1$ such that $\psi_{(t_1,0)}$ has nonreal
  roots, then by continuity, for small enough $t_2$, $\psi_t$ will
  still have nonreal roots with $t=(t_1,t_2)$.  And again, for $t$
  small enough, these zeros will still be distinct.
\end{proof}


\subsection{Some Matrix Inner Functions}
\label{subsec:matrixinner}
We now construct a family of $2\times 2$ matrix valued analytic
functions which are unimodular on the boundary $B$ and have precisely
$6$ zeros in $R$ starting from certain positive $2\times 2$ matrix
valued harmonic functions on $R$.  This will later provide the basis
for our counterexample.

For our purposes, a matrix-valued $F:R\to M_2(\mathbb{C})$ has
\textit{simple zeros} in $R$ if $Z = \{a_1, \dots, a_n\}$ are the
zeros of $\det(F)$, the zeros of $\det(F)$ have multiplicity either
one or two and if the zero $z\in Z$ has multiplicity two, then $F(z) =
0$.

\begin{definition}\rm
 \label{def:standardzero}
 We say that $F$ has a \textit{standard zero set} if $F$ has distinct
 simple zeros $0,a_1,\ldots,a_4\in\mathbb{R}$, with the first having
 multiplicity two and the others multiplicity one.  And if
 furthermore, $\delta_j \neq 0$ are such that $F(a_j)^*\delta_j = 0$,
 $j=1,\ldots, 4$, then no three of the $\delta_j$'s are collinear.
 While it is not required immediately, we also ask that $Ja_j \neq
 P_k$, $k=1,2$, where $P_1,P_2$ are poles of the Fay kernel
 $K^0(\cdot,z)$ (defined in Section \ref{sec:fay-kernels}).  This last
 assumption plays an important role in the diagonalization arguments
 in Section \ref{subsec:diagonalization}.  In this regard it is also
 useful to know that the points $P_1,P_2$ are real and distinct (see
 Lemma \ref{lem:locatecriticalpoints}).  Also not required right away
 but ultimately useful is the assumption that the $\delta_j$'s are
 sufficiently close to the $e_k$'s---how close being determined by
 Theorem \ref{thm:hkernel}.  Here $\{e_1,e_2\}$ denotes the standard
 basis for $\mathbb C^2$, as usual.
\end{definition}

We take $t$ and $p(t)$ as at the end of the last section.

Recall $M(p)$ and $\tau(p)$ from Lemma \ref{lem:allpositive}.  For
notational purposes, write $p(t)=(1,q_1,q_2)$.  Let $\tau$ denote the
common value of $\tau(1,q_1,q_2)$, $\tau(1,q_1^*,q_2)$,
$\tau(1,q_1,q_2^*)$, and $\tau(1,q_1^*,q_2^*)$.

For $0\le \eta\le 1$, let $P^\eta_+$ denote the projection onto the
span of $\eta e_1+(1-\eta^2)^{\frac12}e_2$ and $P^\eta_-$ denote the
projection onto $(1-\eta^2)^{\frac12} e_1 - \eta e_2$ (so $Q^\eta_+
+Q^\eta_- = I$).  Define $H_{\eta,t}$ by
\begin{equation*}
 \begin{split}
   H_{\eta,t}=\tau_0 \mathbb P(\cdot,1)I +& \tau_1 \mathbb
   P(\cdot,c_1-r_1\exp(it_1)) P^1_+
   + \tau_1 \mathbb P (\cdot,c_1-r_1\exp(-it_1)) P^1_-\\
   +&\tau_2\mathbb P(\cdot,-c_2-r_2\exp(it_2)) P^\eta_+ +
   \tau_2\mathbb P(\cdot,-c_2-r_2\exp(it_2)) P^\eta_-.
 \end{split}
\end{equation*}
That is, given $x\in \mathbb C^2$ a unit vector,
$\ip{H_{\eta,t}(\zeta)x}{x}$ is the positive harmonic function on $B$
corresponding to the measure
\begin{equation*}
 \begin{split}
   \mu_{x,x}=\tau_0 \delta_1 +& \tau_1 \left( \delta_{q_1} |x_1|^2 +
     \delta_{q_1^*} |x_2|^2\right) \\
   +&\tau_2\left( \delta_{q_2}|\eta x_1+(1-\eta^2)^\frac12 x_2|^2
     +\delta_{q_2^*} |(1-\eta^2)^\frac12 x_1-\eta x_2|^2\right) .
 \end{split}
\end{equation*}
Then
\begin{equation*}
 \begin{split}
   \int_B Q_j \,d\mu_{x,x} & = \tau_0 Q_j(1) + \tau_1 \left( Q_j(q_1)
     |x_1|^2  +  Q_j(q_1^*)|x_2|^2\right) \\
   & \qquad +\tau_2 \left(Q_j(q_2)|\eta x_1+(1-\eta^2)^\frac12 x_2|^2
     +Q_j(q_2^*) |(1-\eta^2)^\frac12 x_1-\eta x_2|^2\right) .
  \end{split}
\end{equation*}
Using $Q_j(q_k^*)=Q_j(q_k)$ and the choice of $\tau$, gives
\begin{equation*}
 \begin{split}
   \int_B Q_j \,d\mu_{x,x} & = \tau_0 Q_j(1) + \tau_1 Q_j(q_1)
   \left( |x_1|^2  +|x_2|^2\right) \\
   & \qquad +\tau_2 Q_j(q_2)\left( |\eta x_1+(1-\eta^2)^\frac12 x_2|^2
     +
     |(1-\eta^2)^\frac12 x_1-\eta x_2|^2\right) \\
   & = \tau_0 Q_j(1) + \tau_1 Q_j(q_1) + \tau_2 Q_j(q_2) \\
   & = 0.
  \end{split}
\end{equation*}
Hence $\ip{H_{\eta,t}(\zeta)x}{x}$ is the real part of an analytic
function.  It follows that $H_{\eta,t}$ is the real part of an
analytic $2\times 2$ matrix-valued function $G_{\eta,t}(\zeta)$
normalized by $G_{\eta,t}(0)=I$.

Our desired functions are
\begin{equation}
 \label{eq:definePsir}
 \Psi_{\eta,t}=(G_{\eta,t}-I)(G_{\eta,t}+I)^{-1}.
\end{equation}

\begin{lemma}
 \label{lem:PSIONE}
 For each $\eta$,
 \begin{enumerate}
 \item $\Psi_{\eta,t}$ is analytic in a neighborhood of $X$ and
   unitary-valued on $B$;
 \item $\Psi_{\eta,t}(0)=0$;
 \item $\Psi_{\eta,t}(1)=I$;
 \item $\Psi_{\eta,t}(q_1)e_1=e_1$ and $\Psi_{\eta,t}(q_1^*)e_2=e_2$;
 \item $\Psi_{\eta,t}(q_2)(\eta e_1 +(1-\eta^2)^\frac12 e_2)=(\eta e_1
   +(1-\eta^2)^\frac12 e_2)$ and
   $\Psi_{\eta,t}(q_2^*)((1-\eta^2)^\frac12 e_1 - \eta
   e_2)=(1-\eta^2)^\frac12 e_1 - \eta e_2$;
 \item $\Psi_{0,t} =
   \begin{pmatrix}
     \psi_t & 0 \\ 0 & \psi_{-t}
   \end{pmatrix}$ .
 \end{enumerate}
\end{lemma}

\begin{proof}
  In the neighborhood of a point $p\in B$, the Poisson kernel $\mathbb
  P(\zeta,p)$ is the real part of a function of the form
  $g_p(\zeta)(\zeta-p)^{-1}$ where $g_p$ is analytic in a neighborhood
  of $p$ and does not vanish at $p$.  Near any other point $q\in B$,
  $\mathbb P(\zeta,p)$ extends to a harmonic function in a
  neighborhood of $q$ (see Fisher \cite{MR85d:30001}, chapter 4,
  proposition 6.4).  Thus $\mathbb P(\zeta,p)$ is, near $q$, the real
  part of an analytic function, the real part of which is $0$ at $q$.
  
  From the above discussion, if $p\in B$ is different from
  $1,q_1,q_1^*,q_2,q_2^*$, then $G_{\eta,t}$ is, at least locally,
  analytic in a neighborhood of $p$.  Further, $G_{\eta,t}+I$ is
  invertible near $p$, since $G_{\eta,t}(\zeta) =
  H_{\eta,t}(\zeta)+iA(\zeta)$, for some self-adjoint matrix valued
  function $A(\zeta)$ and $H_{\eta,t}(p) = 0$ for such $p$.  Thus,
  $G_{\eta,t}+I$ is invertible at $p$ and by continuity, also
  invertible near $p$.  The relation,
  \begin{equation*}
    I-\Psi_{\eta,t}(\zeta)\Psi_{\eta,t}(\zeta)^* =2(G_{\eta,t}(\zeta)+I)^{-1}
    (G_{\eta,t}(\zeta)+G_{\eta,t}(\zeta)^*)(G_{\eta,t}(\zeta)+I)^{*-1}
  \end{equation*}
  now shows that $\Psi_{\eta,t}$ is unitary at $p$.
 
  From the definition of $G_{\eta,t}$, in a neighborhood of $1$, there
  exists analytic functions $g_1$, $g_2$, $h_1$, and $h_2$ so that the
  real parts of $h_j$ are $0$ at $1$, each $g_j$ is different from $0$
  at $1$, and
  \begin{equation*}
    G_{\eta,t}(\zeta)=  \begin{pmatrix} \frac{g_1}{\zeta-1} & h_1 \\ h_2 &
      \frac{g_2}{\zeta-1} \end{pmatrix}
  \end{equation*}
  Thus,
  \begin{equation*}
   \begin{split}
     &(G_{\eta,t}(\zeta)+I)^{-1} = \frac{1}
     {\frac{g_1+\zeta-1}{\zeta-1}\frac{g_2+\zeta-1}{\zeta-1}-h_1h_2}
    \begin{pmatrix}
      \frac{g_2+\zeta-1}{\zeta-1} & -h_1 \\ -h_2 &
      \frac{g_1+\zeta-1}{\zeta-1} \end{pmatrix} \\
    &= \frac{1}{(g_1+\zeta-1)(g_2+\zeta-1)-h_1h_2(\zeta-1)^2}
      \begin{pmatrix} (g_2+\zeta-1)(\zeta-1)
        & -h_1(\zeta-1)^2 \\ -h_2(\zeta-1)^2 & (g_1+\zeta-1)(\zeta-1)
      \end{pmatrix} .
    \end{split}
  \end{equation*}
  Note that the determinant in the denominator is indeed different
  from $0$ near $1$, and as $\zeta \to 1$, it goes to $g_1(1)g_2(1)
  \neq 0$.  Hence $G_{\eta,t}+I$ is in fact invertible.  Next, by
  directly computing $(G_{\eta,t}-I)(G_{\eta,t}+I)^{-1}=
  \Psi_{\eta,t}$ we see $\Psi_{\eta,t}$ is analytic in a neighborhood
  of $1$ and $\Psi_{\eta,t}(1)=I$.
  
  Now move on to $q_1$.  Near this point there exists analytic
  functions $g$, $h_2$, $h_3$, and $h_4$ so that the real parts of the
  $h_j$ are $0$ at $q_1$, the function $g$ does not vanish at $q_1$
  and
  \begin{equation*}
   G_{\eta,t}=\begin{pmatrix} \frac{g}{\zeta-q_1} & h_2 \\ h_3 & h_4
   \end{pmatrix}
  \end{equation*}
  Since $h_4+1$ has real part $1$ at $q_1$, whereas
  $g(\zeta-q_1)^{-1}$ has a pole and $h_2,h_3$ are analytic at $q_1$,
  we see $G_{\eta,t}+I$ is invertible near $q_1$.  Further, by direct
  computation of $\Psi_{\eta,t}=(G_{\eta,t}-I)(G_{\eta,t}+I)^{-1}$ we
  see that $\Psi_{\eta,t}$ is analytic in a neighborhood of $q_1$ and
  \begin{equation*}
     \Psi_{\eta,t}(q_1)=\begin{pmatrix} 1 & 0 \\ 0 &
       \frac{h_4(q_1)-1}{h_4(q_1)+1} \end{pmatrix} 
  \end{equation*}
  The analogous result holds for the point $q_1^*$.  For the points
  $q_2$, $q_2^*$, the same argument prevails by writing all matrices
  with respect to the orthonormal basis $\{\eta e_1
  +(1-\eta^2)^\frac12 e_2,\eta e_2 +(1-\eta^2)^\frac12 e_1\}$ of
  $\mathbb C^2$.
  
  Finally, (6) is easily seen from the definition of $H_{0,t}$ and the
  functions $G_{0,t}$, $\Psi_{0,t}$.
\end{proof}

We are able to obtain some information on the zeros of
$\Psi_{\eta,t}$, at least for some $\eta$

\begin{lemma}
 \label{lem:chooser}
 There exists an $\epsilon > 0$ such that for all $\eta$ nonzero but
 less than $\epsilon$, $\Psi_{\eta,t}$ has a standard zero set.
\end{lemma}

\begin{proof}
  Recall from our choice of $t$, the zeros of $\psi_t$ are distinct
  and, aside from $0$, are not real.  Of course the zeros of
  $\psi_{-t}$ are just the complex conjugates of those of $\psi_t$.
  Thus
  \begin{equation*}
    \Psi_0 =\begin{pmatrix} \psi_t & 0 \\ 0 & \psi_{-t} \end{pmatrix} ,
  \end{equation*}
  $\Psi_0$ has a standard zero set.
 
  By choosing a sequence $\eta_n$ converging to $0$ and considering
  the uniformly bounded sequence $\Psi_n=\Psi_{\eta_n,t}$, there is
  some subsequence, still written as $\Psi_n$, which converges
  uniformly on compact subsets of $R$ to some $\Psi$.  Consequently,
  \begin{equation*}
    G_n=(I+\Psi_n)(I-\Psi_n)^{-1}
  \end{equation*}
  converges uniformly on compact subsets of $R$ to
  \begin{equation*}
    G=(I+\Psi)(I-\Psi)^{-1}.
  \end{equation*}
  
  Now $H_n$, the real part of $G_n$, is harmonic and
  \begin{equation*}
    H_n-H_0= \mathbb P (\cdot,q_2)(\tau_2(\eta_n) P^{\eta_n}_+ -
    \tau_2(0)P^0_+) + \mathbb P (\cdot,q_2^*) (\tau_2(\eta_n)
    P^{\eta_n}_- -\tau_2(0) P^0_-).
  \end{equation*}
  Thus, as both $(P^{\eta_n}_\pm-P^0_\pm)$ and $(\tau_2(\eta_n) -
  \tau_2(0))$ converge to $0$ with $n$ ($\tau$ is continuous by Lemma
  \ref{lem:allpositive}), the $H_n$ converge to $H_0$ uniformly on
  compact sets.  Since also $G(0)=I=G_0(0)$, it follows that the $G_n$
  converge to $G_0$ and hence the $\Psi_n$ converge uniformly to
  $\Psi_0$ on compact sets.
  
  Let $d_n(\zeta) = \det(\Psi_n(\zeta))$.  This is analytic and
  unimodular on the boundary and $d_n$ converges uniformly to $d_0$ on
  compact sets.  From Rouche's theorem (or else Hurwitz's
  theorem---see Conway \cite{MR80c:30003}), the number of zeros of
  $d_n$ is constant, equaling $6$.  Indeed, these zeros vary
  continuously with $n$ and $d_n$ has a double zero at $0$.  If we let
  $0,w_1,w_2$ denote the zeros of $\psi_t$, then the zeros of
  $\psi_{-t}$ are $0,w_1^*,w_2^*$.  Hence the zeros of $\Psi_0$ are
  exactly given by $\Psi_0(0)=0$, $\Psi_0(w_j)e_1=0$, and
  $\Psi_0(w_j^*)e_2=0$, $\ell = 1,2$, and, listed as
  $0,w_1,w_2,w_1^*,w_2^*$, they are distinct.  Thus, for large $n$,
  the zeros $0,a_1^n,a_2^n,a_3^n,a_4^n$ of $\Psi_n$ satisfy this
  assumption too.
 
  Finally, if $\Psi_n(a_1^n)^*\delta_1^n =0$ and $a_1$ is close to
  $w_1$, then
  \begin{equation*}
    \Psi_0(w_1)^*\delta_1^n  = (\Psi_0(w_1)-\Psi_0(a_1^n))\delta_1^n +
    (\Psi_0(a_1^n)-\Psi_n(a_1^n))\delta_1^n.
  \end{equation*}
  For $n$ large both terms on the right hand side are small and thus
  for $n$ large $\delta$ is close to (a scalar multiple of) $e_1$.
  Thus, for $n$ large enough the $\delta_j^n$ satisfy the assumption
  that no three are collinear.  Finally, since the zeros vary
  continuously with and $\eta$ and for $\eta=0$ the zeros, aside from
  $0$, are not real, whereas $P_1,P_2$ are real (and not zero---see
  Lemma \ref{lem:locatecriticalpoints}), for small enough $\eta$,
  neither $P_1$ nor $P_2$ will be a zero of $\Psi_{\eta,t}$ and the
  nonzero zeros of $\Psi_{\eta,t}$ will be distinct.
  
  This completes the proof.
\end{proof}

\subsection{The Period Matrix and Abel-Jacobi Map}
\label{subsec:period-matrix-abel}
Since the harmonic measures $h_j$ vanish on $B_0$, for $j=1,2$, they
reflect across $B_0$ by
\begin{equation*}
 \begin{split}
   h_j(1/z^*)=&-h_j(z) \\
   h_j(Jz)=&-h_j(z)
 \end{split}
\end{equation*}
to a harmonic function on $Y\setminus B_j$.  Note, when $R^\prime$ is
viewed as the reflection of $R$ across $B_0$, $h_j$ is $-1$ on the
reflection of $B_j$.

While a harmonic conjugate $\tilde{h_j}$ of $h_j$ is only locally
defined, the differential
\begin{equation*}
 \frac{d(h_j+\tilde{h_j})}{dz} \, dz
\end{equation*}
is uniquely determined and globally defined.  If $\zeta$ is a point of
$Y$, then
\begin{equation*}
 \label{eq:defineXi}
  \chi(\zeta) = \frac12 \begin{pmatrix} \int_{-1}^\zeta
    \frac{d(h_1+\tilde{h_1})}{dz} \, dz \\
     \int_{-1}^\zeta \frac{d(h_2+\tilde{h_2})}{dz} \, dz \end{pmatrix}
\end{equation*}
depends upon the path of integration from $-1$ to $\zeta$, but only up
to the number of windings of the path around the boundary components
$B_1,B_2$ and the number of crossings of $B_1,B_2$ (in passing from
$R$ to $R^\prime$).  The choice of $-1$ for the base point is fairly
canonical as it is a Weierstrass point for $Y$ (see subsection
\ref{subsec:minimal}).

The multiple valued nature of $\chi:Y\to \mathbb C^2$ is concisely
described by the period matrix and period lattice for $Y$.  Let
\begin{equation*}
  {\mathbf{P}}_{j,\ell} = \frac{1}{i} \int_{B_\ell}
  \frac{d(h_j+\tilde{h_j})}{dz} \, dz = 
  \int_{B_\ell} \frac{\partial h_j}{\partial n} \, ds.
\end{equation*}
Thus, ${\mathbf{P}}_{j,\ell}$ is the period of the harmonic conjugate
of $h_j$ around $B_\ell$.  The $2\times 2$ matrix ${\mathbf{P}}$ has
positive definite real part and is called the period matrix for $R$.
(Our period matrix differs by a factor of $1/i$ from what most call
the period matrix.)

Let $L$ denote the lattice $\mathbb Z^2+i{\mathbf{P}}\mathbb Z^2$.
The Jacobi variety of $Y$ is the quotient $\mbox{Jac}(Y)=\mathbb
C^2/L$.  Let $[z]$ denote the class of $z\in \mathbb C^2$ in
$\mbox{Jac}(Y)$.  The mapping $\chi_0:Y\to \mbox{Jac}(Y)$,
\begin{equation*}
 \chi_0(\zeta)=[\chi(\zeta)]
\end{equation*}
is well defined and known as the Abel-Jacobi map.  However, it will
often be convenient to work with $\chi$, the multiple valued version
of $\chi_0$.

Let $e_1,e_2$ be the usual basis for $\mathbb C^2$.

\begin{proposition}[\cite{MR93a:30047}, p.~92]
 \label{prop:propertiesofAJ}
 The Abel-Jacobi map has the following properties.
 \begin{enumerate}
 \item $\chi_0$ is a one to one conformal map of $Y$ onto its image in
   $\mbox{Jac}(Y)$;
 \item $\chi_0(J\zeta)=-\chi_0(\zeta)^*$;
 \item If $\zeta\in B$, then $-\chi(\zeta)^*=\chi(\zeta)+n$ for some
   $n\in\mathbb Z^2$.
 \end{enumerate}
\end{proposition}

Property (2) depends upon the choice of basepoint $-1\in B_0$.
Property (3) follows from the observation that for $\zeta\in B$,
$h_j(\zeta)$ is either $0$ or $1$.

\subsection{The Theta Function and its Zeros}
\label{subsec:theta-function-zeros}
Details for most of what follows in this subsection can be found in
Ch.~IV of \cite{MR93a:30047} and Ch.~2 of \cite{MR85h:14026}.

Let $L$ be the lattice defined in the last subsection.  The Riemann
theta function associated to $L$ is the entire function on $\mathbb
C^2$ defined by
\begin{equation*}
 \theta(z)=\sum_{n\in \mathbb Z^2} \exp(-\pi \ip{{\mathbf{P}}n}{n}
 +2\pi i\ip{z}{n}).
\end{equation*}
where $\ip{\cdot}{\cdot}$ is the usual inner product on $\mathbb C^2$.

Straightforward manipulations show $\theta(z^*)=\theta(z)^*$ and
$\theta(-z)=\theta(z)$.

The quasi-periodic behavior of $\theta$ with respect to $L$ is given
by
\begin{equation}
  \label{eq:periodtheta}
   \begin{split}
     \theta(z+\ell)=&\theta(z) \\
     \theta(z+i{\mathbf{P}}m)=&\exp(\pi\ip{{\mathbf{P}}m}{m}-2\pi i
     \ip{z}{m}) \theta(z),
   \end{split}
\end{equation}
where $\ell , m \in {\mathbb{Z}}^2$.

Given $e\in \mathbb C^2$ write $e=u+i{\mathbf{P}}v$, define the theta
function with characteristic $[e]$ by
\begin{equation*}
 \theta[e](z)=\theta \left[\begin{matrix} u \\
     v\end{matrix}\right](z)=\exp (\pi i
 (\ip{{\mathbf{P}}v}{v}+2\ip{u-z}{v})\theta(z-e)).
\end{equation*}
The function $\theta[e]:\mathbb C^2 \to \mathbb C$ obeys the period
laws
\begin{equation}
  \label{eq:periodictranslate}
   \begin{split}
     \theta[e](z+\ell)=&\exp(-2\pi\ip{\ell}{v})\theta[e](z) \\
     \theta[e](z+i{\mathbf{P}}m)=& \exp(2\pi \ip{m}{v})
     \exp(\pi\ip{{\mathbf{P}}m}{m}-2\pi i \ip{z}{m}) \theta[e](z)
   \end{split}
\end{equation}

It turns out that almost all meromorphic functions and differentials
on $Y$ can be represented in terms of translates
$\theta[e](\chi(\zeta)):Y \to {\mathbb C}^2$, which despite being
multiple valued has a well-defined zero set, the description of which
is due to Riemann.

\begin{theorem}
  \label{thm:zerostheta}
  There exists a constant vector $\Delta$ (depending upon the choice
  of basepoint) so that for each $e\in \mathbb C^2$, either
  $\theta[e](\chi(\zeta))$ is identically zero, or
  $\theta[e](\chi(\zeta))$ has exactly $2$ zeros $P_1,P_2$, and
       \begin{equation*}
       \chi(P_1)+\chi(P_2)=e-\Delta \quad \text{modulo }L.
       \end{equation*}
\end{theorem}

Here $\Delta$ is known as the vector of Riemann constants.

The set
\begin{equation*}
 \mathcal{N} = \{e\in \mbox{Jac}(Y): \theta(\chi(\zeta)-e) \mbox{ is
   identically zero}\}
\end{equation*}
is a proper closed subset of $\mbox{Jac}(Y)$.

\subsection{The Prime Form}
\label{subsec:prime-form}
The following fact plays an important role in construction of the
prime form and thus multiple valued meromorphic functions on $Y$ with
prescribed poles and zeros.  Given $e\in\mathbb C^2$ such that
$\theta(e)=0$, define
\begin{equation*}
 \mathcal E_e(\zeta,\xi)=\theta(\chi(\zeta)-\chi(\xi)-e).
\end{equation*}

\begin{theorem}[\cite{MR85h:14026}, Ch.~2, Lemma 3.4]
 \label{thm:zerosoftheta}
 If $e\in\mathbb C^2$, $\theta(e)=0$, and $\mathcal E_e$ is not
 identically zero, then there exists $P\in Y$ so that for each $\xi\in
 Y$, $\xi\neq P$, the zeros of $\theta[e+\chi(\xi)](\chi(\zeta))$,
 which coincide with the zeros of $\mathcal E_e(\zeta,\xi)$, are
 precisely $\xi$ and $P$.
\end{theorem}

The following can be found in Mumford \cite{MR86b:14017} (Lemma~1,
p.~3.208---but see also \cite{MR49:569} and \cite{MR97f:46042}).

\begin{theorem}
  There exists $e_*=\frac12 (u+i{\mathbf{P}}v) \in \mathbb C^2$ such
  that $2e_*=0$ modulo $L$, $\ip{u}{v}$ is odd (equal to $1$ modulo
  $2\mathbb Z$), and $\mathcal E_{e_*}$ is not identically zero.
\end{theorem}

An $e_*$ as in this theorem is called a non-singular odd half period
and for the remainder we take $e_*=\frac12(u_*+i{\mathbf{P}}v_*)$ as
fixed.  Note that $\theta(e_*)=0$, as the fact that $\ip{u_*}{v_*}$ is
an odd integer implies $\theta(e_*)=-\theta(-e_*)$.  In fact,
$e_*+e_*^*=2u_* \in \mathbb Z^2$ so that, using the periodicity of
$\theta$, for $z\in\mathbb C^2$,
\begin{equation*}
 \label{eq:oddhalfperiod}
   \theta(z+e_*^*)=\theta(z+e_*^*+2u_*)=\theta(z-e_*),
\end{equation*}
and so $\theta(e_*)=\theta(-e_*)$ as well, meaning that $\theta(e_*)=
0$.

\begin{lemma}
 \label{lem:PinB}
 There exists a $P\in B$ so that for each $\xi \in R$, the multiple
 valued function $\mathcal E_{e_*}(\zeta,\xi):Y\to \mathbb C$,
 $\mathcal E_{e_*}(\zeta,\xi) = \theta(\chi(\zeta)-\chi(\xi)-e_*)$ is
 not identically zero and has zeros at precisely $P$ and $\xi$.
\end{lemma}

\begin{proof}
  From Theorem \ref{thm:zerosoftheta}, there is a $P$ so that either
  $\mathcal E_{e_*}(\zeta,\xi)$ is identically zero, or has zeros $P$
  and $\xi$.  Accordingly, consider the multiple valued function
  $g:Y\to \mathbb C$ defined by
  $g(\zeta)=\theta(\chi(\zeta)-\chi(-1)-e_*)$.  Since $[\chi(-1)]=0$,
  we may assume $\chi(-1)=0$.  Thus $g(\zeta)=\theta(\chi(\zeta)-e_*)$
  and the fact that $e_*$ is non-singular means $g$ is not identically
  zero.  Hence $g$ has zeros $P$ and $-1$.
  
  Observe, as $e_*+e_*^* \in \mathbb Z^2$, there exists an
  $n\in\mathbb Z^2$ such that $e_*^*=-e_*+n$.  Similarly, as
  $\chi(JP)=-\chi(P)^*$ modulo $L$, there exists $a,b\in \mathbb Z^2$
  so that $\chi(JP)=-\chi(P)^*+a+i{\mathbf{P}}b$.  Hence,
  \begin{equation*}
   \begin{split}
     g(JP)^*&=\theta(\chi(JP)-e_*)^*\\
     &=\theta(-\chi(P)^*+a+i{\mathbf{P}}b-e_*)^* \\
     &=(\kappa \theta(-\chi(P)^*-e_*))^*\\
     &=\kappa^* \theta(-\chi(P)-e_*^*)\\
     &=\kappa^*\theta(-\chi(P)+e_*-n)\\
     &=\kappa^*\theta(-\chi(P) +e_*)\\
     &=\kappa^*\theta(\chi(P)-e_*)\\
     &=\kappa^*g(P),
   \end{split}
  \end{equation*}
  where $\kappa=\exp[\pi\ip{{\mathbf{P}}a}{a}-2\pi
  i\ip{-\chi(P)^*-e_*}{n}]$ is nonzero.  Thus, $g(JP)=0$.  It follows
  that $JP=P$ or $JP=-1$, in which case $P=-1$.  Thus, $JP=P$ and so
  $P$ is in $B$.
  
  Since $P$ is in $B$, $\mathcal E_{e_*}(\cdot,\xi)$ is not
  identically zero and the result follows from Theorem
  \ref{thm:zerosoftheta}.
\end{proof}

Let
\begin{equation*}
   \theta_*(z)=\theta[e_*](z)=\exp[\tfrac12\pi
   i(\ip{{\mathbf{P}}v_*}{v_*}+\ip{u_*-z}{v_*})]\theta(z-e_*) . 
\end{equation*}
If $z,w\in Y$ are different from the $P$ of Lemma \ref{lem:PinB}, then
\begin{equation}
  \label{eq:ratiothetas}
 \frac{\theta_*(\chi(\zeta)-\chi(z))}{\theta_*(\chi(\zeta)-\chi(w))} =
 e^{\pi i[\chi(z) - \chi(w)]}
 \frac{\theta(\chi(\zeta)-\chi(z)-e_*)}{\theta(\chi(\zeta)-\chi(w)-e_*)}
\end{equation}
is multiple valued, but its zero/pole structure is well defined: it
has a zero at $z$ and a pole at $w$.  As we shall see, this will play
an important role in defining reproducing kernels on $R$ with respect
to harmonic measure.

\section{The Fay Kernel Functions of $R$}
\label{sec:fay-kernels}

We now introduce the reproducing kernel $K^a$ as found in Fay
\cite{MR49:569}.  The description of these kernels involves the
critical points for the Green's function $g(\cdot,a)$ for $R$ at the
point $a\in R$.  A point $w$ is a critical point if the gradient of
$g(\zeta,a)$ is $0$ at $w$.  It is well known that in a region of
connectivity $n+1$ there are $n$ of these critical points
(\cite{MR51:13206}, p.~133).  Thus, in $R$ there are two.  In the
sequel we will have use of the following fact about the location of
these critical points for the choice $a=0$.

\begin{lemma}
 \label{lem:locatecriticalpoints}
 The critical points of the Green's function $g(\zeta,0)$ are on the
 real axis, one in each of the intervals, $(-1,c_1-r_1)$ and
 $(c_2+r_2,1)$.
\end{lemma}

\begin{proof}
  For notational ease, let $g(\zeta)=g(\zeta,0)$.  From the symmetry
  of the domain, $\frac{\partial g}{\partial y}=0$ on the $x$-axis in
  $R$.  Since $g$ is $0$ at the points $-1$ and $c_1-r_1$, Rolle's
  Theorem implies there is a point $-1<w_1< c_1-r_1$ so that
  $\frac{\partial g}{\partial x}(w_1)=0$.  Thus, the gradient of $g$
  is zero at $w_1$.  Similarly, there is a point $c_2+r_2 < w_2 < 1$
  such that the gradient of $g$ at $w_2$ is also zero.
\end{proof}

For the remainder of the paper, we let $P_1=Jw_1$ and $ P_2=Jw_2$,
where $w_1,w_2$ are the critical points for the Green's function for
$R$ at $0$.  Thus, $P_1,P_2\in R^\prime$.  These are the points
$P_1,P_2$ which appear in the definition of a standard zero set,
Definition \ref{def:standardzero}.

\begin{theorem}
 \label{thm:FaysK}
 There is a reproducing kernel $K^a$ for the Hardy space $\mathbb
 H^2(R,\omega_a)$ of functions analytic in $R$ with boundary values in
 $L^2(\omega_a)$, where $\omega_a$ is harmonic measure for the point
 $a.$ If $z=a$, then $K^a(\zeta,z)=1.$ Otherwise, $K^a(\zeta,z)$ has
 precisely the poles $P_1(a),P_2(a),Jz$, where $JP_1(a)$ and $JP_2(a)$
 are the critical points for the Green's function for $R$ at $a$, and
 three zeros in $Y$, one of which is at $Ja$.
\end{theorem}

Incidentally, the Hardy space in the theorem corresponds to just one
of the two torus parameter family of rank one bundle shifts over $R$
\cite{MR53:1327}.  Indeed, Ball and Clancey \cite{MR97f:46042} give a
theta function representation for all of the corresponding reproducing
kernels, avoiding the use of the Klein prime form which Fay uses in
his formula for the kernel.

The reader who has skipped sections
\ref{subsec:period-matrix-abel}--\ref{subsec:prime-form} may wish to
skip the proof and proceed directly to subsection
\ref{subsec:appsofthetarep}.

\begin{proof}[Proof of Theorem \ref{thm:FaysK}]
  There is an $e\in \mbox{Jac}(Y)$ so that
  \begin{equation}
    \label{eq:FaysK}
    \begin{split}
      &K^a(\zeta,z)\\
      &=\frac{\theta(\chi(\zeta)+\chi(z)^*+e)\theta(\chi(a)+\chi(a)^*+e)
        \theta_*(\chi(a)+\chi(z)^*) \theta_*(\chi(\zeta)+\chi(a)^*)}
      {\theta(\chi(a)+\chi(z)^*+e)\theta(\chi(\zeta)+\chi(a)^*+e)
        \theta_*(\chi(\zeta)+\chi(z)^*) \theta_*(\chi(a)+\chi(a)^*)}.
    \end{split}
  \end{equation}
  (See Fay \cite{MR49:569}, Proposition~6.15 and Ball and Clancey,
  \cite{MR97f:46042}.)  Straightforward computation using the periodic
  nature of $\theta$ checks that the right hand side is invariant
  under $\chi(\zeta)\mapsto \chi(\zeta)+n+iPm$ for $m,n\in\mathbb C^2$
  so that $K^a(\zeta,z)$ is in fact single valued and meromorphic in
  $Y$.  Further, according to Fay, for $z\in X$, $K^a(\zeta,z)$ is
  analytic as a function of $\zeta \in R$.
  
  One readily verifies from \eqref{eq:FaysK} that $K^a(\zeta,a)=1$, so
  assume $z\neq a$.  Now examine the zero/pole structure of the
  portion of the right side of \eqref{eq:FaysK} depending on $\zeta$:
  \begin{equation}
    \label{eq:ratiosthetas2}
    \frac{\theta(\chi(\zeta)+\chi(z)^*+e)\theta_*(\chi(\zeta)+\chi(a)^*)}
    {\theta(\chi(\zeta)+\chi(a)^*+e)\theta_*(\chi(\zeta)+\chi(z)^*)}.
  \end{equation}
  From the comments following \eqref{eq:ratiothetas},
  $\theta_*(\chi(\zeta)+\chi(a)^*)/\theta_*(\chi(\zeta)+\chi(z)^*)$
  has a zero at $Ja$ and pole at $Jz$.  Each of the remaining theta
  functions has two zeros, which we label $Z_1(z), Z_2(z)$ for the top
  term and $P_1(a), P_2(a)$ for the bottom term.  Hence
  \eqref{eq:ratiosthetas2} (and so \eqref{eq:FaysK}) has zeros at
  $P_1(z), P_2(z), Ja$ and poles at $P_1(a), P_2(a), Jz$, the latter
  which are all in $R'\cup B$, since $K^a$ is analytic.
  
  If any of the poles and zeros were to cancel, either $K^a(\zeta,z)$
  is constant or has two poles and two zeros.  Since these kernels are
  linearly independent and $K^a(\zeta,a)=1$, the first possibility
  cannot occur.  In the second case we would be left with two poles in
  $R'\cup B$, which by Proposition \ref{prop:minimalpoles}, is also
  impossible unless $Jz$ cancels with $Ja$ and the other two poles are
  in $B$.  But then $z=a$ and $K^a(\zeta,a)=1$, which we have already
  ruled out.  Hence $K^a(\zeta,a)$ has order three, with zeros and
  poles as claimed.
\end{proof}


\subsection{Application of the Theta Function Representation of $K^a$}
\label{subsec:appsofthetarep}
 
We assume throughout that for $j= 1,\ldots ,4$, $a_j \in B$ are
distinct and $0\neq\delta_j\in{\mathbb{C}}^2$ have the property that
no three are collinear (the properties of a standard zero set).  As
usual $e_1,e_2$ are the standard basis for ${\mathbb{C}}^2$.  We write
$P_1,P_2$ for $P_1(0), P_2(0)$, the poles of $K^0(\zeta,z)$.

Let $\mathcal M_\delta^0$ denote the span of
\begin{equation*}
  \begin{split}
    & \{K^0(\zeta,0)e_1,K^0(\zeta,0)e_2,K^0(\zeta,a_1)\delta_1,
    \dots,K(\zeta,a_4)\delta_4\} \\
    = & \{e_1,e_2,K^0(\zeta,a_1)\delta_1,\dots,K(\zeta,a_4)\delta_4\}.
  \end{split}
\end{equation*}

Recall that given $p\in \Pi$ there exists an $s\in B_1 \times B_2
=\mathbb T^2$ so that, identifying $s$ with $(1,s)$, $\phi_s =
\phi_p(1)^*\phi_p$.  In particular, the zeros of $\phi_p$,
$z_0^s=0,z_1^s,z_2^s$, depend only upon $s$.

\begin{theorem}
 \label{thm:hkernel}
 Let $a^0_1,\dots,a^0_4$ be points in $R$ so that
 $P_1,P_2,J0,Ja^0_1,Ja^0_2,Ja^0_3,Ja^0_4$ are all distinct.  Let
 $\{e_1,e_2\}$ denote the standard basis for $\mathbb C^2$ and let
 $\delta_1^0=\delta_2^0=e_1$, and $\delta_3^0=\delta_4^0=e_2$.  There
 exists an $\epsilon>0$ so that if
 $|a_j^0-a_j|,\|\delta_j^0-\delta_j\|<\epsilon$, and if
 \begin{equation}
  \label{eq:hkernel}
   h(\zeta)=\sum c_j K^0(\zeta,a_j)\delta_j +v
 \end{equation}
 does not have poles at $P_1,P_2$, then $h$ is constant; i.e., each
 $c_j=0$.
 
 Further, if $h\neq 0$, has a representation as in equation
 (\ref{eq:hkernel}), and if there exists $z_1,z_2\in B$ (not assumed
 to be distinct but not both zero) such that
 \begin{equation*}
   h(\zeta)K^0(\zeta,z_k)=\sum c^k_j K^0(\zeta,a_j)\delta_j +v_k,
 \end{equation*}
 then $h$ is constant, $z_1 = a_{j_1}$, and $z_2 = a_{j_2}$, the
 corresponding $\delta_{j_k}$'s may be taken to be equal to $h$, and
 all the other terms are zero.
\end{theorem}

Note that the theorem is really a statement about the meromorphic
functions $K^0(\zeta,a_j)$ on the double $Y$ and so we view $\zeta$ as
a local coordinate on $Y$.  Indeed, by restricting $\zeta$ to be near
either $P_1$ or $P_2$ it may be assumed that all the points
$\zeta,P_1,P_2,Ja_1,\dots,Ja_4$ are in a single chart $U \subset
R^\prime$ ($U$ is an open simply connected subset of $Y$).

With fixed $a\in R$ distinct from $0,JP_1,JP_2$, where $P_1,P_2$, are
the poles of the kernel $K^0(\zeta,a)$, the residue of the pole of
$K^0(\zeta,a)$ at $P_j$ is given by the value of the analytic function
of $\zeta$
 \begin{equation*}
   (\zeta-P_j)K^0(\zeta,a)
 \end{equation*}
 at the point $\zeta=P_j$.  Let $R_j(a)$ denote this residue.

\begin{lemma}
 \label{lem:continuousresidues}
 The residue $R_j(a)$ varies continuously with $a$.
\end{lemma}

\begin{proof}
  Consider the theta function representation for $K^0(\zeta,a)$ from
  Theorem \ref{thm:FaysK}.  The function
  \begin{equation*}
    f(\zeta)= \theta(\chi(\zeta)+\chi(0)^*+e)
  \end{equation*}
  is analytic and single valued in $U$.  Further, $f(\zeta)$ vanishes
  to order one at $P_j$ and thus can be written as
  \begin{equation*}
    f(\zeta)=(\zeta-P_j)f_j(\zeta),
  \end{equation*}
  where $f_j$ is analytic in $U$ and $f_j(P_j)\ne 0$.  Given a set
  $W\subset U$, let $W^*=\{z^*:z\in W\}$.  Choose neighborhoods $V_j$
  and $W$ of $U$ so that $F:V_j\times W^* \to \mathbb C$ by
  \begin{equation*}
    \begin{split}
      F(\zeta,a^*)=& f(\zeta)K^0(\zeta,a) \\
      =&
      \frac{\theta(\chi(\zeta)+\chi(z)^*+e)\theta(\chi(0)+\chi(0)^*+e)
        \theta_*(\chi(0)+\chi(z)^*) \theta_*(\chi(\zeta)+\chi(0)^*)}
      {\theta(\chi(0)+\chi(z)^*+e) \theta_*(\chi(\zeta)+\chi(z)^*)
        \theta_*(\chi(0)+\chi(0)^*)}
 \end{split}
\end{equation*}
is analytic.  Rewriting gives,
\begin{equation*}
 (\zeta-P_j)K^0(\zeta,a)= \frac{F(\zeta,a)}{f_j(\zeta)}.
\end{equation*}
The lemma now follows from the fact that the right hand side is
analytic in $(\zeta,a^*)\in V_j\times W$.
\end{proof}

\begin{proof}[Proof of Theorem \ref{thm:hkernel}.]
  Without loss of generality, we can assume $\epsilon >0$ is small
  enough that the points $P_1,P_2,Ja_1,\dots,Ja_4$ are all distinct.
  Define
  \begin{equation*}
    R(a_1,a_2)=\begin{pmatrix} R_1(a_1) & R_1(a_2) \\ 
      R_2(a_1) & R_2(a_2) \end{pmatrix}
  \end{equation*}
  and
  \begin{equation*}
    R(a_3,a_4)=\begin{pmatrix} R_1(a_3) & R_1(a_4) \\ 
      R_2(a_3) & R_2(a_4) \end{pmatrix}
  \end{equation*}
  where $R_j(a)$ is the residue of $K^0(\zeta,a)$ at $P_j$ as in the
  Lemma \ref{lem:continuousresidues}.
  
  To prove that $R(a_1,a_2)$ is invertible, let
  \begin{equation*}
    c=\begin{pmatrix} c_1 \\ c_2 \end{pmatrix}
  \end{equation*}
  and $f_c=c_1 K^0(\zeta,a_1) + c_2K^0(\zeta,a_2)$.  Note that
  $R(a_1,a_2)c=0$ if and only if $f_c$ does not have poles at either
  $P_1$ or $P_2$.  In this case, if $f_c$ is not constant, then the
  poles of $f_c$ are precisely the points $Ja_1$ and $Ja_2$ which
  gives the usual contradiction, since both of these points are in
  $R'$.  Thus, $f_c$ is constant.  The kernel $K^0(\zeta,0)=1$, so we
  can express this as
  \begin{equation*}
    0 = c_0 K^0(\zeta,0) + c_1 K^0(\zeta,a_1) + c_2K^0(\zeta,a_2) 
  \end{equation*}
  Since $0,a_1,a_2$ are distinct, the functions
  $K^0(\zeta,0),K^0(\zeta,a_1),K^0(\zeta,a_2)$ are linearly
  independent and hence $c_1=c_2=0$.  Summarizing, if $R(a_1,a_2)c=0$,
  then $c=0$.  It follows that $R(a_1,a_2)$ is invertible and by an
  identical argument, $R(a_3,a_4)$ is invertible.
  
  Consider the function $F$ defined for $\delta_j$ near $\delta_j^0$
  by
  \begin{equation*}
    F=\begin{pmatrix} R_1(a_1)\delta_1 & R_1(a_2)\delta_2  & 
      R_1(a_3)\delta_3 & R_1(a_4)\delta_4 \\
      R_2(a_1)\delta_1 & R_2(a_2)\delta_2  & R_2(a_3)\delta_3 & 
      R_2(a_4)\delta_4 \end{pmatrix}.
  \end{equation*}
  Thus, $F$ takes values in $M_4$, the $4\times 4$ matrices, viewed as
  $2\times 4$ matrices with entries from $\mathbb C^2$.  Clearly $F$
  is continuous in $\delta_j$.  By Lemma \ref{lem:continuousresidues}
  it is also continuous in $a_1,\dots,a_4$.  Indeed, from the form of
  $F$ it is jointly continuous in $a_j,\delta_j$.  Since $F$ is
  invertible at $a_j^0, \delta_j^0$, it follows that there is an
  $\epsilon>0$ so that if
  $|a_j^0-a_j|,\|\delta_j^0-\delta_j\|<\epsilon$, then $F$ is
  invertible.
  
  If $a_j$ and $\delta_j$ are chosen such that $F$ is invertible and
  \begin{equation*}
    h(\zeta)=\sum c_j K^0(\zeta,a_j)\delta_j +v
  \end{equation*}
  does not have poles at $P_j$, then
  \begin{equation*}
    \begin{split}
      0= & \begin{pmatrix} \sum c_j R_1(a_j)\delta_j \\
        \sum c_j R_2(a_j)\delta_j \end{pmatrix}\\
      = & F \begin{pmatrix} c_1\\ c_2 \\ c_3 \\ c_4 \end{pmatrix}.
    \end{split}
  \end{equation*}
  Hence $c=0$ and $h$ is constant.
  
  Now suppose $h\neq 0$ and there exist $z_1,z_2\in B$ (not assumed
  distinct but not both zero) such that
  \begin{equation}
    \label{eq:span}
    h(\zeta)K^0(\zeta,z_k)=\sum c^k_j K^0(\zeta,a_j)\delta_j +v_k,
    \qquad k=0,1,2,
  \end{equation}
  where $z_0 = 0$ (and so $K^0(\zeta, z_0) = 1$.  Using $k=1$ we see
  that $P_1,P_2$ are not poles of $h$, since by the assumptions on the
  distinctness of the $P_k$'s and $a_j$'s, the right side has a pole
  of order at most one at each $P_k$, while the left side has a pole
  of order at least one at these points.  Using $k=0$ we see that $h$
  satisfies the hypothesis of the part of the Theorem which has
  already been proved.  Thus $h$ is constant.
  
  The rest of the result now easily follows using the linear
  independence of the kernels.
\end{proof}


\section{Representing nice matrix valued inner functions}

\subsection{Hahn-Banach Separation}
Recall that $\mathbb H(X)$ denotes the set of functions analytic in a
neighborhood of $X$ and $\mathcal R(X)$ the rational functions with
poles off $X$.  Of particular interest is the set $B\mathbb H(X)$
consisting of those $f\in \mathbb H(X)$ with $\|f\|_R\le 1$.

Let $M_2(\mathbb H(X))$ denote the $2\times 2$ matrices with entries
from $\mathbb H(X)$ and similarly define $M_2(\mathcal R(X))$.  For
$f\in \mathbb H(X)$, $f^*$ denotes its pointwise complex conjugate
while for $F\in M_2(\mathbb H(X))$, $F^*$ is the pointwise adjoint.

For $f, g\in \mathbb H(X)$ or $M_2(\mathbb H(X))$, and $h(z,w) = \sum
f(z)g(w)^*$ (with only finitely many terms) we use the convention,
$\sum f(z)g(w)^*(T) = \sum f(T)g(T)^*$ to define $h(T) = h(T,T^*)$.

Let $\mathcal C$ be the cone generated by
\begin{equation*}
  \{H(z)(1-\psi(z)\psi(w)^*)H(w)^*:\psi \in B\mathbb H(X), H\in
  M_2(\mathbb H(X))\}.
\end{equation*}
Obviously we would get the same set if we were to instead assume the
$H$ is a ${\mathbb{C}}^2$-valued function.

\begin{lemma}
  \label{lem:absorb}
  If $F\in M_2(\mathbb H(X))$, then there exists $\rho>0$ such that
  $I-\rho^2 F(z)F(w)^* \in \mathcal C$.
\end{lemma}

\begin{proof}
  First, suppose
  \begin{equation*}
    F = \begin{pmatrix} f\\ g \end{pmatrix}.
  \end{equation*}
  Choose $0 < \tau$ so big that $\frac{f}{\tau}$ and $\frac{g}{\tau}$
  are in $B\mathbb H(X)$.  Then,
  \begin{equation*}
    \begin{split}
      2\tau^2 I &- F(z)F(w)^* = 2\begin{pmatrix} 1 \\ 0 \end{pmatrix}
      (\tau^2-f(z)f(w)^*)
      \begin{pmatrix} 1 & 0 \end{pmatrix} \\
      &+ 2\begin{pmatrix} 0 \\ 1 \end{pmatrix} (\tau_g^2- g(z)g(w)^*)
      \begin{pmatrix} 0 & 1 \end{pmatrix} +\begin{pmatrix} f(z)\\ 
        -g(z) \end{pmatrix} (1-0) \begin{pmatrix} f(w)^* & g(w)^*
      \end{pmatrix}.
    \end{split}
  \end{equation*}
  Thus $\rho=\frac{1}{\sqrt{2}\tau}$ satisfies the conclusion of the
  lemma.
 
  For general $F$ write $F(z)F(w)^* = G(z)G(w)^*+H(z)H(w)^*$, where
  $G$ and $H$ are the first and second columns of $F$ respectively.
  There exists $\rho_G$ and $\rho_H$ so that both $I-\rho_G^2
  G(z)G(w)^*$ and $I-\rho_H^2 H(z)H(w)^*$ are in $\mathcal C$.  With
  $\rho^2=\frac12 \min \{\rho_G^2,\rho_H^2\}$,
  \begin{equation*}
   \begin{split}
     I-\rho^2 F(z)F(w)^* &= I-\rho^2 G(z)G(w)^* -\rho_2 H(z)H(w)^* \\
     &= \frac12 (I-\rho_G^2 G(z)G(w)^*) +\frac12 (I-\rho_H^2 H(z)H(w)^*)\\
     & + \left( \frac12 \rho_G^2 +\frac12 \rho_H^2 -\rho^2 \right).
   \end{split}
  \end{equation*}
  Each term on the right hand side is evidently in $\mathcal C$.  This
  completes the proof.
\end{proof}

Henceforth, for $F\in M_2(\mathbb H(X))$, we set
\begin{equation*}
  \rho_F = \sup \{\rho>0:I-\rho^2 F(z)F(w)^* \in \mathcal C\}.
\end{equation*}

The following Proposition is an application of the Hahn-Banach
theorem.  It is central to our construction.

\begin{theorem}
  \label{thm:possstatz}
  If there exists a function $F:R\to M_2(\mathbb{C})$ which is
  analytic in a neighborhood of $X$ and unitary-valued on $B$ such
  that $\rho_F <1$, then there exists a Hilbert space $\mathcal H$ and
  an operator $T$ such that $T$ has $X$ as a spectral set, but $T$
  does not have a normal $B$-dilation.
\end{theorem}

\begin{proof}
  The proof features a familiar Hahn-Banach separation argument and
  GNS construction.
  
  From the hypothesis, there exists a $\rho<1$ so that $I-\rho^2
  F(z)F(w)^* \notin \mathcal C$.
 
  Let $\mathcal P$ be the vector space of finite sums
  \begin{equation*}
    \sum h_j(z)g_j(w)^*
  \end{equation*}
  where $H_j,G_j$ are $\mathbb C^2$-valued functions analytic in a
  neighborhood of $X$.  Note, $h(z)g(w)^*$ is pointwise a $2\times 2$
  matrix.
 
  The cone $\mathcal C$ is a convex subset of $\mathcal P$ not
  containing $I-\rho^2 F(z)F(w)^*$ and, by Lemma \ref{lem:absorb}, $I$
  is an internal point of $\mathcal C$.  Hence, there exists a
  nonconstant linear functional $\lambda: \mathcal P \to \mathbb C$ so
  that $\lambda \ge 0$ on $\mathcal C$ and $\lambda(I-\rho^2
  F(z)F(w)^*) \le 0$ we have $\lambda(I)>0$, as otherwise, from Lemma
  \ref{lem:absorb}, $\lambda(H(z)H(w)^*)=0$ for all $H$ and hence
  $\lambda=0$ (see, for example, Holmes \cite{MR53:14085}, \S~11.E).
 
  Let $\mathbb H_2(X)$ denote the $\mathbb C^2$-valued functions
  analytic in a neighborhood of $X$.  For $h,g\in \mathbb H_2(X)$,
  define
  \begin{equation*}
    [h,g]=\lambda( h(z)g(w)^*).
  \end{equation*}
  Since the cone $\mathcal C$ contains $h(z)h(w)^*$, the form
  $[\cdot,\cdot]$ is positive semidefinite on $\mathbb H_2(X)$.
 
  Given $f$ analytic in a neighborhood of $X$, consider the mapping
  $M_f:\mathbb H_2(X) \to \mathbb H_2(X)$ defined by multiplication by
  $f$ so that $M_fg=fg$.  Let $C_f$ be the infimum over all positive
  numbers such that $\frac{f}{C_f}$ takes values in the closed unit
  disk.  Obviously, $C_f = \| f \|_R$.  For any $g\in \mathbb C^2$,
  \begin{equation*}
    (C_f^2 g(z))\left( 1-\frac{f(z)}{C_f} \frac{f(w)^*}{C_f}\right)
    (C_fg(w))^* = g(z)(C_f^2-f(z)f(w)^*)g(w)^* \in \mathcal C.
  \end{equation*}
  Thus,
  \begin{equation*}
    \begin{split}
      C_f^2 [g,g] -[M_fg,M_fg]&=
      C_f^2 \lambda(g(z)g(w)^*)- \lambda( f(z)g(z)g(w)^* f(w))\\
      &= \lambda (g(z) \left (C_f^2-f(z)f(w)^*\right ) g(w)^*) \ge 0,
    \end{split}
  \end{equation*}
  as $\lambda$ is nonnegative on $\mathcal C$.  It follows that each
  $M_f$ defines a bounded operator, still denoted by $M_f$, on the
  Hilbert space $\mathcal H$ obtained from ${\mathbb H}_2(X)$ by
  modding out $[\cdot,\cdot]$ null vectors and completing.
  Furthermore, $\|M_f\| = C_f$.  In particular, with $T=M_{\zeta}$,
  where $\zeta(z)=z$ is the identity function, the set $X$ is a
  spectral set for $T$.  Here we are using $f(T)=M_f$.
 
  To see that $T$ does not have a dilation to a normal operator with
  spectrum in $X$, it suffices to show that $F(T)$ is not a
  contraction for the $F$ in the statement of the theorem.  To this
  end, write
  \begin{equation*}
    F=\begin{pmatrix} F_{11} & F_{12} \\ F_{21} & F_{22} \end{pmatrix}
  \end{equation*}
  so that
  \begin{equation*}
    F^t(T)=\begin{pmatrix} M_{F_{11}} & M_{F_{21}} \\ 
      M_{F_{12}} & M_{F_{22}} \end{pmatrix},
  \end{equation*}
  where $F^t$ denotes the pointwise transpose of $F$.  Let $e_1,e_2$
  denote the (class of) constant functions, $e_j(z)=e_j$, as elements
  of $\mathcal H$.  Compute,
  \begin{equation*}
    \begin{split}
      &\ip{F^t(T)\begin{pmatrix} e_1\\e_2\end{pmatrix}}
      {F^t(T)\begin{pmatrix} e_1\\e_2\end{pmatrix}} =
      \ip{\begin{pmatrix}
          \begin{pmatrix} F_{11}\\ F_{21} \end{pmatrix} \\[6pt]
          \begin{pmatrix} F_{12}\\ F_{22} \end{pmatrix} 
        \end{pmatrix}}
      {\begin{pmatrix}
          \begin{pmatrix} F_{11}\\ F_{21} \end{pmatrix} \\[6pt]
          \begin{pmatrix} F_{12}\\ F_{22} \end{pmatrix} 
        \end{pmatrix}} \\
      & = \lambda\left(
        \begin{pmatrix} F_{11}(z) \\ F_{21}(z) \end{pmatrix}
        \begin{pmatrix} F_{11}(w)^* &  F_{21}(w)^* \end{pmatrix}
        + \begin{pmatrix} F_{12}(z) \\ F_{22}(z) \end{pmatrix}
        \begin{pmatrix} F_{12}(w)^* &  F_{22}(w)^* \end{pmatrix}
      \right)
      \\
      & = \lambda( F(z)F(w)^*).
    \end{split}
  \end{equation*}
  On the other hand,
  \begin{equation*}
    \ip{\begin{pmatrix} e_1\\e_2\end{pmatrix}}
    {\begin{pmatrix} e_1\\e_2\end{pmatrix}}=
    \lambda(e_1 e_1^*)+\lambda(e_2 e_2^*)=\lambda(I).
  \end{equation*}
  Combining the last two equalities gives,
  \begin{equation*}
    \begin{split}
      \ip{(I-F^t(T)^* F^t(T))\begin{pmatrix} e_1\\e_2\end{pmatrix}}
      {\begin{pmatrix} e_1\\e_2\end{pmatrix}}
      =&  \lambda(I-F(z)F(w)^*) \\
      = \lambda\left( \frac{1}{\rho^2} -F(z)F(w)^*\right)-& \left(
        \frac{1}{\rho^2}-1\right) \lambda(I) <0.
    \end{split}
  \end{equation*}
  Therefore, $\|F^t(T)\|>1$.  On the other hand, $\|F^t\|_\infty
  =\sup\{\|F^t(z)\|:z\in X\}$ is the same as $\|F\|_\infty$, since the
  norms of a matrix and its transpose are the same.  In fact, as $F$
  is unitary valued on $B$, so is $F^t$ and thus $\|F^t\|_\infty =1$.
  It now follows from Lemma \ref{lem:lem1} that $T$ does not dilate to
  a normal operator with spectrum in $X$.
\end{proof}

\subsection{Matrix Measures}
This subsection is a brief digression from the main line of
development to collect some needed facts about matrix-valued measures.

Given a compact Hausdorff space $X$, an $m\times m$ matrix-valued
measure
\begin{equation*}
 \mu=\begin{pmatrix} \mu_{j,\ell} \end{pmatrix}_{j,\ell=1}^m
\end{equation*}
on $X$ is an $m\times m$ matrix whose entries $\mu_{j,\ell}$ are
complex-valued regular Borel measures on $X$.  The measure $\mu$ is
positive, written $\mu\ge 0$, if, for each continuous function $f:X\to
\mathbb C^m$,
\begin{equation*}
 f=\begin{pmatrix} f_1\\ \vdots \\ f_m \end{pmatrix}
\end{equation*}
we have
\begin{equation*}
  0\le  \int_X f^* \, d\mu\,  f
   =\sum_{j,\ell} \int_X f_j^* f_\ell \,d\mu_{j,\ell}.
\end{equation*}
The positive measure $\mu$ is bounded by $C>0$ if
\begin{equation*}
 CI_m -\begin{pmatrix} \mu_{j,\ell}(X) \end{pmatrix} \ge 0
\end{equation*}
is positive semidefinite, where $I_m$ is the identity $m\times m$
matrix.

\begin{lemma}
 \label{lem:boundedmu}
 The $m\times m$ matrix-valued measure $\mu$ is positive if and only
 if for each Borel set $\omega$ the $m\times m$ matrix
 \begin{equation*}
  \begin{pmatrix} \mu_{j,\ell}(\omega) \end{pmatrix}
 \end{equation*}
 is positive semi-definite.
 
 Further, if there is a $\kappa$ so that each diagonal entry
 $\mu_{j,j}(X)\le \kappa$, then each entry $\mu_{j,\ell}$ of $\mu$ has
 total variation at most $\kappa$.  In particular, if $\mu$ is bounded
 by $C$, then each entry has variation at most $C$.
\end{lemma}

\begin{proof}
  First suppose $\mu$ is positive.  Let $C(X)$ denote the continuous
  complex-valued functions on $X$.  Fix a vector $c\in \mathbb C^m$.
  Given $f\in C(X)$, the function $cf$ is a continuous $\mathbb
  C^m$-valued function.  Thus, $ \Phi_c : C(X)\to \mathbb C$ given by
  \begin{equation*}
   \Phi_c(f)=\sum \int_X c_j^* c_\ell f \,d\mu_{j,\ell}
  \end{equation*}
  is a positive linear functional.  Hence $\sum c_j^* c_\ell
  \mu_{j,\ell}$ is a positive measure on $X$.  If $\omega$ is a Borel
  set, then
  \begin{equation*}
   \ip{\begin{pmatrix} \mu_{j}{\ell}(\omega) \end{pmatrix} c}{c} 
       = \sum  c_j^* c_\ell \,\mu_{j,\ell}(\omega)>0.
  \end{equation*}
  Since $c$ was arbitrary, the matrix in the lemma is positive
  semidefinite for each $\omega$.
 
  Conversely, suppose the matrix in the lemma is positive semidefinite
  for each Borel set $\omega$.  If $f:X\to \mathbb C^m$ is a
  measurable simple functions, $f=\sum v_j \chi_{\omega_j}$, then
  \begin{equation*}
    \int_X f^* \, d \mu \, f =\sum_j \ip{\mu(\omega_j)v_j}{v_j} \ge 0.
  \end{equation*}
  For more general $f:X\to \mathbb C^m$, choose a sequence of
  measurable simple functions converging to $f$ pointwise.  Then,
  \begin{equation*}
   0\le  \int_X f_n^* \, d \mu \, f_n \rightarrow 
   \int_X f^* \, d\mu \, f.
  \end{equation*}
  
  Now suppose that for each diagonal entry $\mu_{j,j}(X)\le \kappa$.
  For $g\in C(X)$, let $\mu_{j,\ell}(g)$ denote the integral of $g$
  with respect to $d\mu_{j,\ell}$ and for $j\ne \ell$, let
  \begin{equation*}
    \lambda=-\frac{\mu_{j,\ell}(g)^*}{|\mu_{j,\ell}(g)^*|} 
      \|g\|_\infty.
  \end{equation*}
  Observe
  \begin{equation*} 
    \begin{split}
      0\le & \ip{\begin{pmatrix} \mu_{j,j} & \mu_{j,\ell}\\
          \mu_{\ell,j} & \mu_{\ell,\ell}\end{pmatrix}
        \begin{pmatrix} \lambda \\ g \end{pmatrix}}
      {\begin{pmatrix} \lambda \\ g \end{pmatrix}} \\
      =& \|g\|_\infty^2 \mu_{j,j}(X)+\mu_{\ell,\ell}(|g|^2)
      - 2 |\mu_{j,\ell}(g)|\|g\|_\infty \\
      \le & 2(\|g\|_\infty^2 \kappa -|\mu_{j,\ell}(g)|\|g\|_\infty).
    \end{split}
  \end{equation*}
  It follows that $\mu_{j,\ell}$ is a continuous linear functional on
  $C(X)$ with norm at most $\kappa$.  In particular, the variation of
  $\mu_{j,\ell}$ is at most $\kappa$.
 
  Finally note that by choosing $c=e_j$, where $e_j$ is the $j$-th
  standard basis vector for $\mathbb C^m$, it follows that each
  $\mu_{j,j}$ is a positive measure.  Further, if we now suppose $\mu$
  is bounded by $C$, then $\mu_{j,j}(X)\le C$.
\end{proof}

\begin{lemma}
 \label{lem:convergemu}
 If $\mu^n$ is a sequence of positive $m\times m$ matrix-valued
 measures on $X$ which are all bounded above by $C$, then there is a a
 positive $m\times m$ matrix-valued measure $\mu$ on $X$ also bounded
 above by $C$.  Hence there is a subsequence $\mu^{n_k}$ of $\mu^n$
 converging to $\mu$ weak-$*$; i.e., for each pair of continuous
 functions $f,g:X\to \mathbb C^m$,
 \begin{equation*}
  \sum_{j,\ell} \int_X f_\ell g_j^* \,d\mu^n_{j,\ell} \to 
     \sum_{j,\ell} \int_X f_\ell g_j^* \,d\mu_{j,\ell}.
 \end{equation*}
\end{lemma}

\begin{proof}
  By the previous lemma, for each $j,\ell$ the measures
  $\mu^n_{j,\ell}$ are bounded in variation by $C$.  Hence, we can
  find a subsequence, which for convenience, we will still denote by
  $\mu^n$ so that each $\mu^n_{j,\ell}$ converges weak-$*$ to some
  $\mu_{j,\ell}$ with variation at most $C$.
  
  Suppose $f:X\to \mathbb C^m$ is continuous.  We have
  \begin{equation*}
    0 \le \int_X f^* \, d\mu^n f
    = \sum_{j,\ell} \int_X f_j^* f_\ell \, d\mu_{j,\ell}^n
    \to \sum_{j,\ell} \int_X f_j^* f_\ell \, d\mu_{j,\ell}
    = \int_X f^* \, d\mu \, f.
  \end{equation*}
  Hence $\mu$ is a positive measure.  Further, for a vector $c\in
  \mathbb C^m$ thought of as a constant function,
  \begin{equation*}
      0 \le C\|c\|^2 - \ip{\left( \mu_{j,\ell}^n(X) \right) c}{c}
      \to C\|c\|^2 - \ip{\left( \mu_{j,\ell}(X)\right) c}{c} .
  \end{equation*}
  Thus, $\mu$ is bounded above by $C$.
\end{proof}

\begin{lemma}
 \label{lem:rnmu}
 If $\mu$ is a positive $m\times m$ matrix-valued measure on $X$, then
 the diagonal entries, $\mu_{j,j}$ are positive measures.  Further,
 with $\nu=\sum_{j=1}^m \mu_{j,j}$, there exists an $m\times m$
 matrix-valued function $\Delta:X\to M_m(\mathbb C)$ so that
 $\Delta(x)$ is positive semidefinite for each $x\in X$ and $d\mu=
 \Delta \, d\nu$; i.e., for each pair of continuous functions
 $f,g:X\to \mathbb C^m$,
  \begin{equation*}
   \sum_{j,\ell} \int_X  g_j^* f_\ell \,d\mu _{j,\ell} =
      \sum_{j,\ell} \int_X g_j^* \Delta_{j,\ell}\, f_\ell  \,d\nu.
  \end{equation*}
\end{lemma}

\begin{proof}
  From Lemma \ref{lem:boundedmu}, if $\omega$ is a Borel set and
  $\mu_{j,j}(\omega)=0$, then $\mu_{j,\ell}(\omega)=0$ for each $j$.
  Thus, each $j,\ell$ the measure $\mu_{j,\ell}$ is absolutely
  continuous with respect to $\nu$.  By the Radon-Nikodym Theorem,
  there exists $\nu$ integrable functions $\Delta_{j,\ell}$ so that
  $d\mu_{j,\ell}=\Delta_{j,\ell} \, d\nu$.
  
  Once again fix a vector $c\in\mathbb C^m$.  By Lemma
  \ref{lem:boundedmu}, for each Borel set $\omega$,
  \begin{equation*}
    0\le \ip{\left( \mu_{j,\ell}(\omega) \right) c}{c}
    = \ip{\left( \int_\omega \Delta_{j,\ell} \, d\nu \right) c}{c}
    = \int_\omega \sum_{j,\ell} c_j^* c_\ell \Delta_{j,\ell} \,d\nu.
  \end{equation*}
  Thus, $\sum_{j,\ell} c_j^* c_\ell \Delta_{j,\ell} \ge 0$ almost
  everywhere with respect to $\nu$.
 
  Choose a countable dense subset $\{c^n\}$ of $\mathbb C^m$.  For
  each $n$ there is a set $E_n$ of $\nu$-measure zero such that off of
  $E_n$ the function $\sum_{j,\ell} (c_j^n)^* c_\ell^n
  \Delta_{j,\ell}$ is non-negative.  For $x\in X \setminus (\cup
  E_n)$, we have $\sum_{j,\ell} (c_j^n)^* c_\ell^n \Delta_{j,\ell}(x)
  \ge 0$ for each $n$.  By continuity of the inner product in $\mathbb
  C^m$, it follows that $\sum_{j,\ell} c_j^* c_\ell \Delta_{j,\ell}(x)
  \ge 0$ for all $c\in\mathbb C^m$ for almost all $x$; that is, the
  matrices $\Delta(x) =(\Delta_{j,\ell}(x))$, for $x\in X \setminus
  (\cup E_n)$, are positive semidefinite.
\end{proof}

\subsection{Representations in terms of the $\phi_p$ }
The following is a companion to Proposition \ref{thm:possstatz}.
Given a subset $S\subset R$, a function $\Gamma:S\times S\times\Pi \to
\mathbb C$ is a positive kernel if for each $p\in\Pi$, the matrix
\begin{equation}
  \label{eq:poskernel}
  (\Gamma(z,w;p))_{z,w\in S'},
\end{equation}
$S'$ a finite subset of $S$.  If $S = R$, then $\Gamma$ is analytic if
$\Gamma(z,w;p)$ is analytic in $z$ and conjugate analytic in $w$ for
all $p$.

The point of the next result is that we do not know \textit{a priori}
that $1-FF^* \in \mathcal C$ even if $\rho_F = 1$.

\begin{proposition}
  \label{prop:tightrep}
  Suppose $F$ is a $2\times 2$ matrix-valued function analytic in a
  neighborhood of $R$, $F$ is unitary-valued on $B$, and $F(0) = 0$.
  If $\rho_F = 1$ and if $S\subset R$ is a finite set, then there
  exists a probability measure $\mu$ on $\Pi$ and a positive kernel
  $\Gamma:S\times S\times \Pi\to \mathbb C$ so that
  \begin{equation*}
    1-F(z)F(w)^* = \int_{\Pi} (1-\phi_p(z)\phi_p(w)^*)\Gamma(z,w;p)
    \,d\mu(p)
  \end{equation*}
  for all $z,w\in S$.
\end{proposition}

\begin{proof}
  Choose a sequence $0 < \rho_n < 1$ such that $\rho_n$ converges to
  $1$.  For each $n$, there exists vector functions $H_{n,j}$ and
  functions $\psi_{n,j}$ analytic in a neighborhood of $R$ such that
  $\psi_{n,j}\in B \mathbb H(X)$ and
  \begin{equation*}
    1-\rho_n^2 F(z)F(w)^* = \sum_{j = 1}^{N_n}
    H_{n,j}(z)(1-\psi_{n,j}(z)\psi_{n,j}(w)^*)H_{n,j}(w)^*. 
  \end{equation*}
  By post composition with a M\"obius transformation if necessary, it
  may be assumed without loss of generality that $\psi_{n,j}(0) = 0$
  for each $n,j$.
  
  For each $n,j$, there exists $h_{n,j}(z,p)$, analytic as a function
  of $z$ in a neighborhood of $R$, and a probability measure
  $\nu_{n,j}$ on $\Pi$ so that by Proposition \ref{eq:scalarrep1},
  \begin{equation*}
    1-\psi_{n,j}(z)\psi_{n,j}(w)^* = \int_{\Pi} h_{n,j}(z,p)
    (1-\phi_p(z)\phi_p(w)^*) h_{n,j}(w,p)^* \,d\nu_{n,j}(p).
  \end{equation*}
  Observe, as all the $\psi$ and $\phi$ vanish at $0$,
  \begin{equation*}
    1 = \int h_{n,j}(0,p) h_{n,j}(0,p)^* \,d\nu_{n,j}(p) .
  \end{equation*}
  Let
  \begin{equation*}
    \nu_n = \sum_{j=1}^{N_n} \nu_{n,j}.
  \end{equation*}
  By the Radon-Nikodym theorem, there exists a nonnegative function
  $u_{n,j}(p)$ such that
  \begin{equation*}
    d\nu_{j,n} = u_{n,j}(p)^2 \,d\nu_n.
  \end{equation*}
  Thus
  \begin{equation*}
    \begin{split}
      1-\psi_{n,j}(z)&\psi_{n,j}(w)^*\\
      &= \int_{\Pi} u_{n,j}(p)h_{n,j}(z,p) [1-\phi_p(z)\phi_p(w)^*]
      h_{n,j}(w,p)^*u_{n,j}(p) \,d\nu_{n}(p).
    \end{split}
  \end{equation*}
  
  Let
  \begin{equation*}
    \Gamma_n(z,w;p)= \sum_j H_{j,n}(z)\,u_{n,j}(p)\,h_{n,j}(z,p)\,
    h_{n,j}(w,p)^* \,u_{n,j}(p)^*\,H_{n,j}(w)^*.
  \end{equation*}
  By construction, $\Gamma_n$ is analytic in $z$ conjugate analytic in
  $w$ in a neighborhood of $R$, is positive semidefinite as a kernel,
  and
  \begin{equation}
   \label{eq:tightrep1}
    I-\rho_n^2 F(z)F(w)^* = \int_{\Pi} [1-\phi_p(z)\phi_p(w)^*]
    \Gamma_n(z,w;p)\,d\nu_n(p).
  \end{equation}
  
  For fixed $z$, $I-\rho^2_n F(z)F(z)^* \le I$ and there exists a
  $\epsilon_z$ such that $1-|\phi_p(z)|^2 \ge \epsilon_z$.  Thus,
  \begin{equation*}
    I\ge \epsilon_z \int_{\Pi} \Gamma_n(z,z;p)\,d\nu_n(p).
  \end{equation*}
  
  Let $C$ denote the maximum of the set $\{\epsilon_z^{-1}:z\in S\}$
  and $m$ the cardinality of $S$.  The sequence of $m\times m$
  measures with $2\times 2$ entries
  \begin{equation*}
    d\mu_n = {\left( \Gamma_n(z,w;p)\, d\nu_n(p) \right) }_{z,w\in S}
  \end{equation*}
  are positive and the diagonal entries are bounded by $C$.  A
  positive $k\times k$ matrix whose diagonal entries are at most $C$
  is bounded above by $kCI_k$.  Thus, it follows from the results of
  the previous section that there exists a positive measure $\nu$ and
  a pointwise positive definite matrix valued function
  \begin{equation*}
   \Gamma(p)={\left( \Gamma(z,w;p) \right)}_{z,w\in S}
  \end{equation*}
  so that some subsequence of $\mu_n$ converges to $\Gamma \, d\nu$
  weak-$*$, where $\nu$ can be taken to be a probability measure by
  scaling $\Gamma$ if necessary.  For notational ease, we continue to
  denote the subsequence by $\mu_n$.
  
  For $z,w\in S$ fixed, the expression $(1-\phi_p(z)\phi_p(w)^*)$ is
  continuous in $p$ by Lemma \ref{lem:extremeclosed}.  Thus, letting
  $n$ tend to infinity in equation (\ref{eq:tightrep1}) gives
  \begin{equation*}
    I-F(z)F(w)^* = \int_{\Pi} [(1-\phi_p(z)\phi_p(w)^* ]
    \Gamma (z,w;p)\,d\nu (p).
   \end{equation*}
\end{proof}

\subsection{Transfer Function Representations}

For present purposes, a unitary colligation $\Sigma = (U,K,\mu)$
consists of a probability measure $\mu$ on $\Pi$, a Hilbert space $K$,
and a unitary a unitary operator $U$ on the direct sum
$(L^2(\mu)\otimes K)\oplus \mathbb C^2$, written as
\begin{equation}
  \label{eq:matrixU}
  U = \begin{pmatrix} A & B \\ C & D \end{pmatrix}
\end{equation}
with respect to the direct sum decomposition.  Here $L^2(\mu) \otimes
K$ signifies $K$-valued $L^2(\mu)$.

Define $\Phi:R\to B(L^2(\mu)\otimes K)$ by $(\Phi(z) f)(p) =
\phi_p(z)f(p)$.  Of course $\Phi$ depends upon $\mu$ and $K$, but this
dependence is suppressed.  The transfer function associated to
$\Sigma$ is
\begin{equation}
  \label{eq:transfer}
  W = W_\Sigma (z) = D+C\Phi(z)(I-A\Phi(z))^{-1} B.
\end{equation}
Note that as $A$ must be a contraction and $\Phi(z)$ is a strict
contraction, the inverse in \eqref{eq:transfer} exists for $z\in R$.
Moreover, since $\Phi(z)(I-A\Phi(z))^{-1} = (I-\Phi(z)A)^{-1}\Phi(z)$,
the transfer function of $\Sigma$ may also be expressed as
\begin{equation}
  \label{eq:transferalt}
  W = D+C(I-\Phi(z)A)^{-1} \Phi(z)B.
\end{equation}

\begin{proposition}
  \label{prop:transferrep}
  The transfer function is contractive-valued, $\|W_\Sigma(z)\|\le 1$
  for all $z\in R$.  Indeed, for $z,w \in R$,
  \begin{equation*}
    I-W_\Sigma(z) W_\Sigma(w)^* = C(I-\Phi(z)A)^{-1}
    (I-\Phi(z)\Phi(w)^*) (I-\Phi(w)A)^{*-1}C^*.
  \end{equation*}
\end{proposition}

By now the proof is entirely standard.  Simply use the equation
(\ref{eq:transfer}), and $DB^* = -CA^*$, $BB^* = I-AA^*$, and $DD^* =
I-CC^*$ to verify:
\begin{eqnarray*}
  I-W(z)W(w)^*& = & 1-[D+C\Phi(z) (I-A\Phi(z) )^{-1} B]
  [D^*+B^*(I-\Phi^*(w)A)^{-1}\Phi^*(w) C] \\
  &=& 1- DD^* - [C\Phi(z) (I-A\Phi(z) )^{-1} B] [B^*(I-\Phi^*(w) A^*
  )^{-1}\Phi^*(w) C]\\ 
   &&\qquad - DB^*(I-\Phi^*(w) A^* )^{-1}\Phi^*(w) C^* - C\Phi(z)
   (I-A\Phi(z) )^{-1} BD^* \\ 
  &=& C\,[ 1 - \Phi(z) (I-A\Phi(z) )^{-1}(1-AA^*)(I-\Phi^*(w) A^*
  )^{-1}\Phi^*(w)\\ 
   &&\qquad + A^*(I-\Phi^*(w) A^* )^{-1}\Phi^*(w) + \Phi(z)
   (I-A\Phi(z) )^{-1} A ]\,C^* \\ 
  &=& C\,[ 1 - (I-\Phi(z) A )^{-1}\Phi(z)(1-AA^*)\Phi^*(w)(I-
  A^*\Phi^*(w) )^{-1} \\ 
   &&\qquad + A^*\Phi^*(w) (I- A^*\Phi^*(w) )^{-1} + (I-\Phi(z) A
   )^{-1}\Phi(z) A ]\,C^* \\ 
  &=& C(I-\Phi(z) A )^{-1}\,[ (I-\Phi(z) A )(I- A^*\Phi^*(w) ) -
  \Phi(z)(1-AA^*)\Phi^*(w)\\ 
   &&\qquad + (I-\Phi(z) A )A^*\Phi^*(w) + \Phi(z) A(I- A^*\Phi^*(w) ) ]
   \,(I- A^*\Phi^*(w) )^{-1}C^* \\
  &=& C(I-\Phi(z) A )^{-1}[ 1-\Phi(z)\Phi^*(w) ](I- A^*\Phi^*(w) )^{-1}C^*.
\end{eqnarray*}

Note that pointwise on $R$, we can define $H(w) =
(I-A^*\Phi(w)^*)^{-1}C^* :\mathbb C^2 \to L^2(\mu)\otimes K$.  Thus,
for $w$ fixed, $H(w)^*$ is a function on $\Pi$, which we emphasize by
writing as $H_p(w)^*$.  Interpreting the representation in Proposition
\ref{prop:transferrep} in terms of the space $L^2(\mu)\otimes K$
gives,
\begin{equation*}
  I-W(z)W(w)^* = \int (1-\phi_p(z)\phi_p(w)^*)H_p(z)H_p(w)^* \,d\mu(p).
\end{equation*}

\subsection{Nevanlinna-Pick Interpolation}

The following proposition is an Agler-Pick type interpolation theorem
for some matrix-valued functions on $R$.  The proof proceeds via a
transfer function realization for a solution.  This is a rip-off of
methods pioneered by Agler, followed by Ball and others, and now
standard.  We eventually show that, roughly speaking, this theorem
only applies to $2\times 2$ matrix functions which are, up to a fixed
unitary, the direct sum of scalar contractive functions.

\begin{proposition}
  \label{prop:NPinterpolate}
  If $S\subset R$ is a finite set, $W:S\to M_2(\mathbb{C})$, and if
  there is a positive kernel $\Gamma:S\times S \times \Pi \to
  M_2(\mathbb{C})$ such that
  \begin{equation*}
    I-W(z)W(w)^* = \int_{\Pi} (1-\phi_p(z)\phi_p(w)^*)\Gamma(z,w;p)
    \,d\mu(p)
  \end{equation*}
  for all $z,w \in S$, then there exists $G:R\to M_2(\mathbb{C})$ such
  that $\|G(z)\|\le 1$ and $G(z) = W(z)$ for $z\in S$.  Indeed, there
  exists a finite dimensional Hilbert space $K$ (dimension at most
  twice the cardinality of $S$) and a unitary colligation $\Sigma =
  (U,K,\mu)$ so that
  \begin{equation*}
    G = W_\Sigma,
  \end{equation*}
  and hence there exists $\Delta:R\times R\times \Pi \to
  M_2(\mathbb{C})$ a positive analytic kernel such that
  \begin{equation*}
    I-G(z)G(w)^* = \int_{\Pi} [1-\phi_p(z)\phi_p(w)^*]\Delta (z,w;p)
    \,d\mu(p)
  \end{equation*}
  for all $z,w\in R$.
\end{proposition}

\begin{proof}
  Once again, this is by now standard.  For $p\in \Pi$, the rank of
  the block matrix with $2\times 2$ matrix entries
  \begin{equation*}
    (\Gamma(z,w;p))_{z,w\in S}
  \end{equation*}
  is at most $2N$, where $N$ is the cardinality of $S$.  Thus, by
  Kolmogorov's theorem (see, for example, \cite{MR2003b:47001},
  Theorem 2.53, especially the second proof), there exists a Hilbert
  space $K$ of dimension $2N$ and a function $H:S\to L^2(\mu)\otimes
  B(\mathbb C^2,K)$, denoted $H_p(z)$, such that $\Gamma(z,w;p) =
  H_p(z)H_p(w)^*$ ($\mu$ almost everywhere).
  
  Let $\mathcal E$ and $\mathcal F$ denote the subspaces of
  $(L^2(\mu)\otimes K) \oplus \mathbb C^2$ spanned by
  \begin{equation*}
    \left\{ \begin{pmatrix} H_s(w)^* x \\ W(w)^*x \end{pmatrix}: x\in
    \mathbb C^2, w\in S\right\},
  \end{equation*}
  and
  \begin{equation*}
    \left\{ \begin{pmatrix} \phi_s(w)^* H_s(w)^* x \\ x \end{pmatrix}:
      x\in \mathbb C^2, w\in S \right\}
  \end{equation*}
  respectively.  The mapping $V$ from $\mathcal E$ to $\mathcal F$
  determined by
  \begin{equation*}
    V \begin{pmatrix} H_s(w)^* x \\ W(w)^*x \end{pmatrix} =
    \begin{pmatrix} \phi_s(w)^* H_s(w)^* x \\ x \end{pmatrix}
  \end{equation*}
  is an isometry since
  \begin{equation*}
    \left< \begin{pmatrix} H_s(w)^* x \\ W(w)^*x \end{pmatrix},
    \begin{pmatrix} H_s(z)^* y \\ W(z)^*y \end{pmatrix}\right> 
     = \left< \int H_s(z)H_s(w)^*\, d\mu(s)\, x,y\right> +
    \left< W(z)W(w)^*x,y \right>
  \end{equation*}
  and
   \begin{equation*}
    \begin{split}
      &\left< \begin{pmatrix} \phi_s(w)^* H_s(w)^* x \\ x
        \end{pmatrix},
       \begin{pmatrix} \phi_s(z)^* H_s(z)^* y \\ y \end{pmatrix}\right> \\
     =& \left< \int \phi_s(z)\phi_s(w)^* H_s(z) H_s(w)^*\,d\mu(s)\,
       x,y\right> + \left< x,y\right> .
   \end{split}
  \end{equation*}
  Both $\mathcal E$ and $\mathcal F$ are finite dimensional, so there
  exists a unitary $U$ on $(L^2(\mu)\otimes K )\oplus \mathbb C^2$
  such that $U^*$ extends $V$.  Let $\Sigma = (U,K,\mu)$ denote the
  resulting unitary colligation.
   
  Write $U$ as in equation (\ref{eq:matrixU}).  Since $U^*$ restricted
  to $\mathcal E$ is $V$,
  \begin{equation*}
    \begin{pmatrix} A^* & C^*\\ B^* & D^*\end{pmatrix}
    \begin{pmatrix} \phi_s(w)^* H_s(w)^* x \\ x \end{pmatrix} 
     = \begin{pmatrix} H_s(w)^* x \\ W(w)^* x\end{pmatrix}.
  \end{equation*}
  Expressing this as a system of equations
  \begin{equation*}
    \begin{split}
      A^* \phi_s(w)^*H_s(w)^*x + C^*x & = H_s(w)^* x \\
      B^* \phi_s(w)^*H_s(w)^*x + D^*x & = W(w)^* x.
    \end{split}
  \end{equation*}
  Solving the first equation for $H_s(w)^*x$ gives,
  \begin{equation*}
    H_s(w)^*x = (I-A^*\Phi(w)^*)^{-1} C^* x.
  \end{equation*}
  Substituting into this the second equation now gives,
  \begin{equation*}
    B^*\Phi(w)^* (I-A^*\Phi(w)^*)^{-1} C^* x = W(w)^*x.
  \end{equation*}
  It follows that for each $z\in S$,
  \begin{equation*}
    W_\Sigma(z) = W(z).
  \end{equation*}
\end{proof}

Next is a uniqueness result for Nevanlinna-Pick interpolation on $R$.

\begin{proposition}
  \label{prop:unique}
  Suppose $F:R\to M_2(\mathbb{C})$ is analytic in a neighborhood of
  $X$, unitary on $B$, and with a standard zero set.  Then there
  exists $S\subset R$ a set with seven elements such that, if $Z:R\to
  M_2(\mathbb{C})$ is contractive-valued and $Z(z) = F(z)$ for $z\in
  S$, then $Z = F$.
\end{proposition}

\begin{proof}
  Let $K^0$ denote the Fay kernel for $R$ defiend in Section
  \ref{sec:fay-kernels}.  That is, $K^0$ is the reproducing kernel for
  the Hilbert space ${\mathbb H}^2(R)$ of functions analytic in $R$
  with $L^2$ boundary values with respect to harmonic measure on $B$
  with respect to the point $0$.  Let ${\mathbb H}^2_2(R)$ denote
  $\mathbb C^2$-valued ${\mathbb H}^2(R)$.  Since $F$ is
  unitary-valued on $B$, the mapping $V$ on ${\mathbb H}_2^2(R)$
  defined by $VG(z) = F(z)G(z)$ is an isometry.  As is shown below,
  the kernel of $V^*$ is the span of $\{K^0(\cdot, a_j)\gamma_j :j =
  1, 2, \dots, 6\}$ where $F(a_j)^* \gamma_j = 0$ and, of course,
  $\gamma_j \neq 0$; that is, the pair $(a_j, \gamma_j)$ is a zero of
  $F^*$.
  
  Before proceeding, we note that if $\varphi$ is scalar valued and
  analytic in a neighborhood of $R$, has no zeros on $B,$ and has
  distinct zeros $w_1,\dots,w_n\in R$ of multiplicity one, and if
  $f\in {\mathbb H}^2(R)$ with $f(w_j)=0$, then $f=\varphi g$ for a
  $g\in {\mathbb H}^2(R)$.
  
  Now suppose $\psi\in {\mathbb H}^2(R)$ and for all $h\in {\mathbb
    H}^2(R)$, we have $\ip{\psi}{\varphi h}=0$.  Then there is a
  linear combination $f=\psi -\sum_1^n c_j K^0(\cdot,w_j)$ so that
  $f(w_j)=0$ as the set $\{K^0(\cdot,w_j):1\le j \le n\}$ is linearly
  independent, and thus by the above remark, $f=\varphi g$ for some
  $g\in {\mathbb H}^2(R)$.  Since $\ip{K^0(\cdot,w_j)}{\varphi h} =
  \varphi(w_j)^*h(w_j)^* =0$ for each $j$ and $h$, it follows that
  $\ip{f}{\varphi h}=0$ for all $h$.  Choosing $h = g$ gives,
  $\ip{\varphi g}{\varphi g}=0$ from which it follows that $g=0$.
  Hence $\psi$ is in the span of $\{K^0(\cdot,w_j):1\le j\le n\}$.
  This shows $\{K^0(\cdot,w_j):1\le j\le n\}$ is a basis for the
  orthogonal complement of $\{\varphi h:h\in {\mathbb H}^2(R)\}$.
  
  We next determine the kernel of $V^*$.  Write $a_5=a_6=0$.  Since
  $F(0)=0$, there is a function $H$ analytic in a neighborhood of $X$
  so that $F=zH$.  The function $\varphi=z\det(H)$ satisfies the
  hypothesis of the preceding paragraph.
  
  Let
  \begin{equation*}
    G = \begin{pmatrix} h_{22} & -h_{12} \\
      -h_{21} & h_{11} \end{pmatrix},
  \end{equation*}
  where $H=(h_{j,\ell})$.  Verify $FG=zHG= z\det(H)I$, where $I$ is
  the $2\times 2$ identity matrix.
  
  Now suppose $x\in {\mathbb H}_2^2(R)$ and $V^*x = 0$.  Let $x_1,
  x_2$ denote the coordinates of $x$.  For each $g\in {\mathbb
    H}^2_2(R)$,
  \begin{equation*}
     0=\ip{Gg}{V^*x}=\ip{VGg}{x}
     =\ip{z\det(H)g}{x}
     = \ip{z\det(H)g_1}{x_1} +  \ip{z\det(H)g_2}{x_2}.
  \end{equation*}
  It follows from the discussion above that each $x_j$ is in the span
  of $\{K^0(\zeta,a_j):1\le j\le 6\}$ and therefore $x$ is in the span
  of $\{K^0(\zeta,a_j)v:1\le j\le 6, v\in \mathbb C^2\}$.  In
  particular, there exists vectors $v_j \in \mathbb C^2$ such that
  \begin{equation*}
    x = \sum_1^6  v_j K^0(\cdot,a_j).
  \end{equation*}
  Since, as is readily verified, $V^* vK^0(\cdot,a) = F(a)^*v
  K^0(\cdot,a)$ and $F(0)^* = 0$,
  \begin{equation*}
    0 = V^* x = \sum_1^4 F(a_j)^* v_j K^0(\cdot,a_j).
  \end{equation*}
  But $K^0(\cdot,a_j)$, $j=1\ldots 4$ are linearly independent, and so
  $F(a_j)^*v_j = 0$ for each $j$.  Conversely, if $F(a_j)^*v_j=0$,
  then $V^* v_j K^0(\cdot,a_j)=0$ so that $\{K^0(\zeta,a_j)v_j:1\le
  j\le 6\}$ is a basis for the kernel of~$V^*$.
  
  The projection onto the kernel of $V^*$ is $I-VV^*$ and since the
  kernel of $V^*$ has dimension six, $I-VV^*$ has rank six.  Thus, for
  any finite set $A\subset R$, the block matrix with $2\times 2$
  matrix entries
  \begin{equation*}
   \begin{split}
     M_A &= {\left({\left(\left< (I-VV^*)K^0(\cdot,w)e_j,
               K^0(\cdot,z)e_\ell\right>
           \right)}_{j,\ell = 1,2}\right)}_{z,w\in A}\\
     &= ((I-F(z)F(w)^*)K^0(z,w))_{z,w\in A}
   \end{split}
  \end{equation*}
  has rank at most six.  (Here $\{e_1, e_2\}$ is the usual basis for
  $\mathbb C^2$.)  In particular, if $A = \{a_1, \dots,a_6\}$ where
  the $a_j$ are the zeros of $\det(F)$ (with $a_5=a_6=0$), then $M_A$
  has exactly rank six.  Choose points $a_7,a_8$ distinct from
  $a_1,\dots,a_6$ so that the set $S=\{a_1,\dots,a_6,a_7,a_8\}$ has
  exactly seven distinct points.  Since $S$ contains $A$, the rank of
  $M_S$ is at least six.  On the other hand, by what is proved above,
  it has rank at most six.  Thus the rank of $M_S$ is six.
  
  Since the matrix $M_S$ is $14\times 14$ (viewed as a $7\times 7$
  block matrix with $2\times 2$ matrix entries), and $M_S$ has rank
  six, the kernel of $M_S$ has dimension eight.  Further, as the
  dimension of the subspace $\mathcal L_1 =\{\gamma\otimes
  e_1:\gamma\in \mathbb C^7\}$,
  \begin{equation*}
    \mathcal L_1=\left \{ \begin{pmatrix} 
      \begin{pmatrix} \gamma_1 \\ 0 \end{pmatrix} \\[6pt]
        \begin{pmatrix} \gamma_2 \\ 0 \end{pmatrix} \\
        \vdots \\
     \begin{pmatrix} \gamma_7 \\ 0 \end{pmatrix} \end{pmatrix}:
      \gamma=\begin{pmatrix} \gamma_1\\ \gamma_2\\ \vdots\\ \gamma_7 
           \end{pmatrix} \in\mathbb C^7
         \right \}.
  \end{equation*}
  of $\mathbb C^7 \otimes \mathbb C^2$ has dimension seven, it follows
  that the there exists a nonzero $x_1=y_1\otimes e_1$ in $\mathcal L$
  and in the kernel of $M_S$.  Similarly, by considering $\mathcal L_2
  = \{\gamma\otimes e_2:\gamma\in \mathbb C^7\}$, there exists a
  nonzero $x_2$ in the kernel of $M_S$ of the form $x_2 = y_2\otimes
  e_2$.
  
  Let $X = \begin{pmatrix} x_1 & x_2 \end{pmatrix}$.  Thus, $X$ is a
  $14\times 2$ matrix,
  \begin{equation*}
    X= \begin{pmatrix} 
      \begin{pmatrix} (y_1)_1 & 0\\ 0 &  (y_2)_1 \end{pmatrix} \\[6pt]
        \begin{pmatrix} (y_1)_2 & 0 \\ 0& (y_2)_2 \end{pmatrix} \\
        \vdots \\
     \begin{pmatrix} (y_1)_7 & 0 \\ 0 & (y_2)_7 \end{pmatrix} 
   \end{pmatrix}.
  \end{equation*}
  It is convenient to use $S$ to index itself, so that we have
  \begin{equation*}
    X(w) = \begin{pmatrix} x_1(w) & x_2(w) \end{pmatrix}
     = \begin{pmatrix} y_1(w) & 0 \\ 0 & y_2(w) \end{pmatrix}
  \end{equation*}
  for the $w\in S$ coordinate of $X$.  The identity $M_SX = 0$
  becomes,
  \begin{equation*}
    \sum_{w\in S} K^0(z,w)X(w) = F(z) \sum_{w\in S} K^0(z,w)F(w)^* X(w)
  \end{equation*}
  for each $z\in S$.
  
  Now suppose $Z:R \to M_2(\mathbb{C})$ is analytic, contractive
  valued, and $Z(z) = F(z)$ for $z\in S$.  The operator $W$ of
  multiplication by $Z$ on ${\mathbb H}^2_2(R)$ is a contraction and
  \begin{equation*}
    W^* K^0(\cdot,w)x = Z(w)^* x K^0(\cdot,w).
  \end{equation*}
  Given $\zeta \in R$, $\zeta \notin S$, let $S^\prime = S\cup
  \{\zeta\}$ and consider the decomposition of
  \begin{equation*}
    N_{\zeta} =
    {\left( {\left(\ip{(I-Z(z)Z(w)^*)K^0(\cdot,w)e_j}{K^0(\cdot,z)e_\ell}
          \right) }_{j,\ell}\right)}_{z,w\in S^\prime}
  \end{equation*}
  into blocks according to $S$ and $\{\zeta \}$.  Thus, $N_\zeta$ is
  an $8\times 8$ matrix with $2\times 2$ block entries.  The upper
  left $7\times 7$ block matrix, the block determined by $S,$ is $M_S$
  since $Z(z) = F(z)$ for $z\in S$.  Let
  \begin{equation*}
    Y = \begin{pmatrix}  X \\ \begin{pmatrix} 0&0\\0 &0 \end{pmatrix}
    \end{pmatrix}.
  \end{equation*}
  Since $N_\zeta$ is positive semi-definite and $M_S X = 0$, it
  follows that $N_\zeta Y = 0$.  An examination of the last two
  entries (the last $2\times 2$ block) of the product $N_\zeta Y = 0$
  gives,
  \begin{equation}
    \label{eq:solveF}
  \sum_{w\in S}
    K^0(\zeta,w)X(w) = Z(\zeta) \sum_{w\in S} Z(w)^* K^0(\zeta,w)X(w).
  \end{equation}
  
  From the form of $X$ and since the kernel functions
  $\{K^0(\cdot,w):w \in S\}$ form a linearly independent set in
  ${\mathbb H}^2(R,\omega_0)$, it follows that the right hand side in
  (\ref{eq:solveF}) has rank two for all but at most countably many
  $\zeta$.  Hence,
  \begin{equation*}
    \sum_{w\in S} K^0(z,w)F(w)^* X(w)
  \end{equation*}
  also has rank two and $Z$ is completely determined by the identity
  (\ref{eq:solveF}).  Since $F$ is also determined by this identity,
  $Z = F$.
\end{proof}

The fact that we used $K^0$ in the last proof is not significant: any
reproducing kernel would have worked.  However, certain facts that
came to light in the proof will play a role in the proof of Lemma
\ref{lem:rank}, where Fay kernels are needed.

\subsection{A tight representation for some $F$}

The results in the previous sections now combine to produce a tight
representation for $F:R\to M_2(\mathbb{C})$ which is analytic across
the boundary, unitary on the boundary, has at most six simple zeros in
$R$, and has $\rho_F =1$.

\begin{theorem}
  \label{thm:rep}
  Suppose $F$ is a $2\times 2$ matrix-valued function analytic in a
  neighborhood of $R$, $F$ is unitary-valued on $B$, with a standard
  zero set.  If $\rho_F = 1$, then there exists a unitary colligation
  $\Sigma = (U,K,\mu)$ such that $F = W_\Sigma$ and so that the
  dimension of $K$ is at most $14$.  In particular, $\mu$ is a
  probability measure $\mu$ on $\Pi$, and there is an analytic
  function $H:R\to L^2(\mu)\otimes M_{14,2}(\mathbb{C})$, denoted
  $H_p(z)$ so that
  \begin{equation*}
    1-F(z)F(w)^* = \int_{\Pi} (1-\phi_p(z)\phi_p(w)^*) H_p(z)H_p(w)^*
    \,d\mu(p)
  \end{equation*}
  for all $z,w\in R$.
\end{theorem}

\begin{proof}
  Choose, using Proposition \ref{prop:unique}, a finite set $S$ of $R$
  with the property: if $G:R\to M_2(\mathbb{C})$ is analytic,
  contractive-valued, and $G(z) = F(z)$ for $z\in S$, then $G = F$.
  
  Using this $S$, Proposition \ref{prop:tightrep} produces a
  probability measure $\mu$ and positive kernel $\Gamma: S\times
  S\times \Pi \to M_2(\mathbb{C})$ such that
  \begin{equation*}
    1-F(z)F(w)^* = \int_{\Pi} (1-\phi_p(z)\phi_p(w)^*)\Gamma(z,w;p)
    \,d\mu(p)
  \end{equation*}
  for all $z,w\in S$.
  
  By Proposition \ref{prop:NPinterpolate} there exists a unitary
  colligation $\Sigma = (U,K,\mu)$ so that $W_\Sigma (z) = F(z)$ for
  $z\in S$.  From our choice of $S$, we see that $F = W_\Sigma$.  The
  integral representation follows.
\end{proof}

\subsection{Fay Kernels Reprise}
Suppose $F$ satisfies the hypothesis of Theorem \ref{thm:rep} and let
$\mu$ denote the measure which appears in the conclusion.

Recall Fay's variant of the Szeg{\H o} kernel for $R$, $K^0(\zeta,z)$.
Also recall, for each $p\in \Pi$ there is an $s\in \{1\}\times
B_1\times B_2$ so that $\phi_s(z)\phi_s(w)^*=\phi_p(z)\phi_p(w)^*$.
In this case we label the zeros of $\phi_p$ (which of course are the
same as those of $\phi_s$) $z_0^s(=0),z_1^s,z_2^s$.

\begin{lemma}
  \label{lem:rank}
  Suppose $F$ satisfies the hypothesis of Theorem \ref{thm:rep}: $F$
  is a $2\times 2$ matrix-valued function analytic in a neighborhood
  of $R$, $F$ is unitary-valued on $B$, with a standard zero set, and
  $\rho_F = 1$.  Assume furthermore that it is represented as in the
  conclusion of Theorem \ref{thm:rep}.  Let $a_5=a_6=0$ and
  $\delta_5=e_1$ and $\delta_6=e_2$.  Then there exists a set $E$ of
  $\mu$ measure zero such that for $p\notin E$, for each $v\in
  {\mathbb C}^{14}$ and for $\ell = 0,1,2$, the vector function
  $H_p(\zeta)v K^0(\zeta,z_\ell)$ is in the span of
  $\{K^0(\zeta,a_j)\delta_j\}$, where $z_0(=0),z_1,z_2$ are the zeros
  of $\phi_p$.  As a consequence, $H_p$ is analytic on $R$ and extends
  to a meromorphic function on $Y$.
\end{lemma}

\begin{proof}
  Given a finite set $Q\subset R$,
  \begin{equation*}
    M_Q = {\left( (I-F(z)F(w)^*)K^0(z,w)\right)}_{z,w\in Q}
  \end{equation*}
  has rank at most six.  Moreover, as interpreted in the proof of
  Proposition \ref{prop:unique}, the range of $M_Q$ lies in the span
  of $\{ (K^0(z,a_j)\delta_j )_{z\in Q}:j = 1,\dots,6\}$.  Here
  $(K^0(z,a_j)\delta_j )_{z\in Q}$ is a column vector indexed by $Q$.
  
  From the representation for $F$ from Theorem \ref{thm:rep},
  \begin{equation*}
    M_Q = {\left(\int H_p(z)(1-\phi_p(z)\phi_p(w)^*)K^0(z,w) H_p(w)^*
        \,d\mu(p) \right)}_{z,w\in Q}.
  \end{equation*}
  
  For each $p$, multiplication by $\phi_p$, denoted $M_p$, is
  isometric on $\mathbb H^2(R,\omega_0)$, the Hilbert space of
  analytic functions on $R$.  Hence $1-M_pM_p^* \geq 0$, as is
  $(I-M_pM_p^*)\otimes E$, where $E$ is the $n\times n$ matrix
  consisting of all ones.  From the reproducing property of
  $K^0(\cdot,z)$, $M_p^* K^0(\cdot,z)=\phi_p(z)^* K^0(\cdot,z)$.
  Thus, if $Q$ is a set of $n$ points in $R$ and $x$ is the vector
  $(K^0(\cdot,w))_{w\in Q}$, then
  \begin{equation*}
    P_Q(p) = \ip{(M_p\otimes E)x}{x} 
    ={\left((1-\phi_p(z)\phi_p(w)^*)K^0(z,w)\right)}_{z,w\in Q} \geq
    0.
  \end{equation*}
  
  If we set $\tilde Q = Q \cup \{z_j(s)\}$, for any $j=0$, $1$, or
  $2$, then $P_{\tilde Q}(p) \geq 0$ as well.  Furthermore, the upper
  left $n\times n$ block equals $P_Q(p)$ and the right $n\times 1$
  column is $(K_0(z,z_j(s)))_{z\in Q}$.  Hence as a vector
  $(K_0(z,z_j(s)))_{z\in Q} \in \ran P_Q(p)^{1/2} = \ran P_Q(p)$ for
  $j=0,1,2$.
  
  Since $P_Q(p) \geq 0$,
  \begin{equation*}
    N_Q(p) = {\left( H_p(z)(1-\phi_p(z)\phi_p(w)^*)K^0(z,w) H_p(w)^*
      \right)}_{z,w\in Q}
  \end{equation*}
  is also positive semidefinite for each $p$.  If $M_Q x = 0$, then
  \begin{equation*}
    0 = \int \left< N_Q(p)x,x\right> \,d\mu(p),
  \end{equation*}
  so that $\left< N_Q(p)x,x\right> = 0$ for almost every $p$.  Since
  $N_Q(p)$ is positive semidefinite, $N_Q(p)x = 0$ almost everywhere.
  Choosing a basis for the kernel of $M_Q$, it follows that there is a
  set $E_Q$ of measure zero so that for $p\notin E_Q$, the kernel of
  $M_Q$ is a subspace of the kernel of $N_Q(p)$.  Thus, for such $p$,
  the range of $N_Q(p)$ is a subspace of the range of $M_Q$.  In
  particular, the rank of $N_Q(p)$ is at most $6$.
  
  Moreover, if we let $D_{Q}(p)$ denote the ($2\times 14$ block)
  diagonal matrix with $(z,z)$ entry $H_p(z)$ ($z\in Q$), then $N_Q(p)
  = D_Q(p)P_Q(p)D_Q(p)^*$.  Since $P_Q(p)$ is positive semidefinite,
  we conclude that the range of $D_Q(p)P_Q(p)$ is a subspace of $M_Q$.
  Thus, as $(K^0(z,z_j(s)))_{z\in Q}$ is in the range of $P_Q(p)$,
  $(H_p(z)v K^0(z,z_j(s)))_{z\in Q}$ is in the range of $M_Q$ for
  every $v\in \mathbb C^{14}$ and $j=0,1,2$.
  
  Now suppose $Q_n \subset R$ is a finite set with $Q_n \subset
  Q_{n+1}$, $Q_0 = \{a_1,\ldots ,a_4,a_5(=0)\}$, and $\mathcal D =
  \bigcup_n Q_n$ a determining set.  Then since
  \begin{equation}
    (H_p(z)v K^0(z,z_j(s)))_{z\in Q_n} \in \ran M_{Q_n} \subseteq
    \bigvee_j (K^0(z,a_j)\delta_j )_{z\in Q},
  \end{equation}
  there are constants $c_j^n(p)$ such that
  \begin{equation}
    \label{eq:restrictH}
    H_p(z)v K^0(z,z_j(s)) = \sum_{j=1}^5 c_j^n(p) K^0(z,a_j)\delta_j,
    \qquad z\in Q_n.
  \end{equation}
  By linear independence of the $K^0(\cdot,a_j)$'s, the $c_j^n(p)$'s
  are uniquely determined when $n=0,1,\ldots$ by \eqref{eq:restrictH}.
  Since $Q_{n+1} \supset Q_n$, we must in this case have $c_j^{n+1}(p)
  = c_j^n(p)$ for all $n$, and thus there are unique constants
  $c_j(p)$ such that
  \begin{equation}
    H_p(z)v K^0(z,z_j(s)) = \sum_{j=1}^5 c_j(p) K^0(z,a_j)\delta_j,
    \qquad z\in \mathcal D ,\qquad p=0,1,2.
  \end{equation}
  By considering this equation with $j=5$, and using the fact that
  $K^0(z,0) = 1$, we see that $H_p$ agrees with an analytic function
  on the determining set $\mathcal D$.  It follows that we can assume
  that $H_p(z)$ is analytic for each $p\notin E$ and that the relation
  of equation (\ref{eq:restrictH}) holds throughout $R$.  Furthermore,
  since the $K^0(\cdot,a_j)$'s extend to meromorphic functions on $Y$,
  $H_p$ does as well.
\end{proof}

\subsection{Diagonalization}
\label{subsec:diagonalization}
In this subsection we show that when $\rho_F = 1$, a contractive
matrix valued function with unitary boundary values and a standard
zero set is diagonalizable.

\begin{lemma}
  Suppose $F$ is a $2\times 2$ matrix valued function on $R$ whose
  determinant is not identically zero.  If there exists a $2\times 2$
  unitary matrix $U$ and scalar valued functions $\psi_1,\psi_2:R\to
  \mathbb C$ such that $F(z)F(w)^*=U D(z)D(w)^*U^*$, where
 \begin{equation*}
   D=\begin{pmatrix} \psi_1 & 0\\ 0&\psi_2 \end{pmatrix},
 \end{equation*}
 then there exists a unitary matrix $V$ so that $F=UDV$.
\end{lemma}

\begin{proof}
  The hypothesis imply $D(z)^{-1}U^*F(z)=D(w)^*U^* F(w)^{*-1}$
  whenever $F(z)$ and $F(w)$ are invertible.  Hence,
  $D(z)^{-1}U^*F(z)=V$ is constant.  One readily verifies that
  $V^*V=I$.
\end{proof}

\begin{theorem}
 \label{thm:diagonal}
 Suppose $F$ is a $2\times 2$ matrix valued function which is analytic
 in a neighborhood of $R$, unitary valued on $B$ and has a standard
 zero set, $\delta_j,a_j$ with the following property.  If $h$
 satisfies
 \begin{equation*}
   h(\zeta) =\sum_1^4 c_j K^0(\zeta,a_j)\delta_j +v
 \end{equation*}
 for some $c_1,\dots,c_4\in\mathbb C$ and $v\in \mathbb C^2$ and if
 $h$ does not have pole at either $P_1$ or $P_2$, then $h$ is
 constant.
 
 If $\rho_F = 1$, then $F$ is diagonalizable; i.e., there exists a
 $2\times 2$ unitary matrices $U$ and $V$ and analytic functions
 $\varphi_j:R\to \mathbb C$ such that, with
 \begin{equation*}
   D=\begin{pmatrix} \varphi_1 & 0\\ 0 & \varphi_2 \end{pmatrix},
 \end{equation*}
 $F=VDU$.
\end{theorem}

\begin{proof}
  By Lemma \ref{lem:rank}, we may assume that, except perhaps on a set
  $\Delta_0$ of measure $0$, if $h$ is a column of some $H_p$, then
  $h(\zeta)K^0(\zeta,z_\ell(s))$ is in the span of
  $\{K^0(\zeta,a_j)\delta_j:1\le j\le 6\}$ for $\ell=0,1,2$.  Here
  $z_0(p)=0,z_1(p),z_2(p)$ are the zeros of $\phi_p$.  By hypothesis,
  $h$ (and so $H_p$) is constant as a function of $\zeta$.  From our
  normalization of the domain (see Corollary \ref{cor:notriplezero}),
  one of the zeros of $\phi_p$, say $z_1(p)$, is not zero.  Thus,
  Theorem \ref{thm:hkernel} (or rather the argument there) also
  implies, if $h$ is not zero, then there is a $1\le j_1(p) \le 4$ so
  that $z_1(p)=a_{j_1(p)}$ and $h$ is a multiple of $\delta_{j_1(p)}$.
  Thus, every column of $H_p$ is a multiple of $\delta_{j_1(p)}$.
  
  Theorem \ref{thm:rep} gives us the representation
  \begin{equation}
    \label{eq:con1}
    1-F(z)F(w)^* = \int_{\Pi} (1-\phi_p(z)\phi_p(w)^*)
    H_pH_p^*\,d\mu(p).
  \end{equation}
  Substituting $w=0$ we find
  \begin{equation}
    \label{eq:con2}
    I=\int H_p H_p^*\,d\mu(p).
  \end{equation}
  In view of equation (\ref{eq:con2}), rearranging equation
  (\ref{eq:con1}) gives,
  \begin{equation}
    \label{eq:con3}
    F(z)F(w)^* =\int_\Pi \phi_p(z)\phi_p(w)^* H_p H_p^* d\mu(p).
  \end{equation}
  
  Since the columns of $H_p$ are all multiples of the single vector
  $\delta_{j_1(p)}$, $H_pH_p^*$ is rank one and thus may be written as
  $G(p)G(p)^*$, for a vector $G(p)\in\mathbb C^2$.  (Indeed, $G(p)$ is
  the square root of the sum of the squares of the norms of the
  columns of $H_p$ times $\delta_{j_1(p)}$.) Consequently,
  \begin{equation}
    \label{eq:con4}
    F(z)F(w)^* =\int_{\Pi} \phi_p(z)\phi_p(w)^* G(p)G(p)^*\, d\mu(p).
  \end{equation}
  Since $F(a_j)^*\delta_j=0$ for all $j$, equation (\ref{eq:con4})
  gives,
  \begin{equation*}
    0=\delta_j^* F(a_j)F(a_j^*)\delta_j 
   =\int |\phi_p(a_j)|^2 \|G(p)^* \delta_j\|^2 \, d\mu(p).
  \end{equation*}
  Consequently, for each $j$, $\phi_p(a_j)^*G(p)^* \delta_j=0$ for
  almost every $p$, and so off of a set $Z_0 \subset \Pi$ of measure
  zero, $\phi_p(a_j)^*G(p)^* \delta_j=0$ for all $p$ and each $j$.
  Thus, by defining $G(p)=0$ on $Z_0$, we may assume that equation
  (\ref{eq:con4}) holds and
  \begin{equation*}
    \phi_p(a_j)^*G(p)^* \delta_j=0
  \end{equation*}
  for all $p,j$.
  
  Let $\Delta_0=\{p\in \Pi: G(p)=0\}$.  If $p\notin \Delta_0$, then
  for each $j$, either $\phi_p(a_j)=0$ or $G_p^* \delta_j=0$.  Since
  $G_p$ is a multiple of $\delta_{j_1(p)}$ and no three of the
  $\delta_j$ are collinear, it follows that $\phi_p$ has zeros at two
  of the $a_j$, say $a_{j_1(p)}, a_{j_2(p)}$, and the
  $\delta_{j_3(p)}$ and $\delta_{j_4(p)}$ are collinear (and
  orthogonal to $\delta_{j(p)}$, where
  $\{a_1,\dots,a_4\}=\{a_{j_1(p)},\dots, a_{j_4(p)}\}$.  We can now
  return to Theorem \ref{thm:hkernel} and conclude that
  $z_1(p)=a_{j_1(p)}, z_2(p)=a_{j_2(p)}$.  In particular, $\phi_p$ has
  distinct zeros and $\delta_{j_1(p)}$ and $\delta_{j_2(p)}$ are
  collinear (and orthogonal to $\delta_{j_3(p)}$ and
  $\delta_{j_4(p)}$).  Let $J_1=\{a_{j_1(p)},a_{j_2(p)}\}$ and
  $J_2=\{a_{j_3(p)},a_{j_4(p)}\}$.  In addition, let $\Delta_1$ denote
  the one dimensional subspace of $\mathbb C^2$ spanned by
  $\delta_{j_1(p)}$ and $\Delta_2$ the the one dimensional subspace of
  $\mathbb C^2$ spanned by $\delta_{j_3(p)}$.
 
  If $p^\prime \notin \Pi_0$, then by arguing as above, either
  $G(p^\prime)\in \Delta_1$ or $G(p^\prime)\in \Delta_2$.  In the
  former case, the zeros of $\phi_{p^\prime}$ are in $J_2$ and in the
  later in $J_1$.  Hence, for each $p$, either
  \begin{itemize}
  \item[(0)] $G(p)=0$; or
  \item[(1)] $G(p)\in \Delta_1$ and the nonzero zeros $z_1(p),z_2(p)$
    are in $J_2$; or
  \item[(2)] $G(p)\in \Delta_2$ and the nonzero zeros $z_1(p),z_2(p)$
    are in $J_1$.
  \end{itemize}
  Let
  \begin{equation*}
    \begin{split}
      \Pi_0=&\{p\in \Pi: \mbox{ (0) holds }\}, \\
      \Pi_1 = &\{p\in\Pi: \mbox{ (1) holds }\}, \\
      \Pi_2 =& \{p\in\Pi: \mbox{ (2) holds }\}.
    \end{split}
  \end{equation*}
  If $p,q\in \Pi_1$, then $\phi_p$ and $\phi_q$ have the same zeros,
  and are therefore equal up to a rotation.  Hence, for $p,q\in
  \Pi_1$, $\phi_p(z)\phi_p(w)^*=\phi_q(z)\phi_q(w)^*$.  Choose $p^1\in
  \Pi_1$ and let Let $\psi_1=\phi_{p^1}$ denote a representative.  If
  $\Pi_2$ is not empty, choose $p^2\in \Pi_2$ and let
  $\psi_2=\phi_{p^2}$.  Otherwise, let $\psi_2=0$.  Substituting into
  equation (\ref{eq:con3}) and writing the integral as the sum of the
  integrals over $\Pi_1$ and $\Pi_2$ gives,
  \begin{equation*}
    F(z)F(w)^* = h_1 \psi_1(z)\psi_1(w)^* h_1 + h_2
    \psi_2(z)\psi_2(w)^* h_2^*
  \end{equation*}
  for some $h_j\in \Delta_j$.  Substituting $z=w=1$ and using the fact
  that $F(1)F(1)^*=I$ shows $\{h_1,h_2\}$ is an orthonormal basis for
  $\mathbb C^2$ (and $\psi_2\ne 0$).  Thus, we can apply the previous
  lemma and conclude that $F$ is diagonalizable.
\end{proof}

\section{The Obstruction}
We now demonstrate a $2\times 2$ matrix function $F$ unitary on the
boundary, analytic across the boundary, and with a standard zero set
which cannot be diagonalized.  Since by Theorem \ref{thm:diagonal},
any function satisfying these conditions which has $\rho_F=1$ is
diagonalizable, $\rho_F$ must be less than $1$ for this $F$ and so by
Theorem \ref{thm:possstatz}, rational dilation does not hold.

Recall the matrix inner functions $\Psi_{\eta,t}$ introduced in
subsection \ref{subsec:matrixinner}.  The $t$ was fixed at the outset
of that section.  By Lemma \ref{lem:chooser} for small $\eta$,
$\Psi=\Psi_{\eta,t}$ has a standard zero set.  Let $a_j^0,\delta_j^0$
denote the standard zero set for $\Psi_{0,t}$.  (In particular, we can
assume $\delta_1^0=\delta_2^0=e_1$ and $\delta_3^0=\delta_4^0=e_2$.)
Let $\epsilon^0 >0$ denote the $\epsilon>0$ in Theorem
\ref{thm:hkernel} corresponding to this zero set.  Thus, for a small
$\eta>0$ the function $\Psi=\Psi_{\eta,t}$ has a standard zero set and
this zero set satisfies the conditions in Theorem \ref{thm:hkernel}.
This $\Psi$ is our funny function.

\begin{lemma}
 \label{lem:propsoffunnypsi}
 $\Psi$ has the following properties:
 \begin{enumerate}
 \item $\Psi$ is unitary valued on $B$;
 \item $\Psi$ has a standard zero set;
 \item the zero set $a_1,\dots,a_4$, $\delta_1,\dots,\delta_4$ has the
   property: If $h$ satisfies
   \begin{equation*}
     h(\zeta) =\sum_1^4 c_j K^0(\zeta,a_j)\delta_j +v
   \end{equation*}
   for some $c_1,\dots,c_4\in\mathbb C$ and $v\in \mathbb C^2$ and if
   $h$ does not have pole at either $P_1$ or $P_2$, then $h$ is
   constant.
 \item $\Psi$ is not diagonalizable; i.e., there does not exist fixed
   unitaries $U$ and $V$ so that $U\Psi V^*$ is pointwise diagonal.
 \end{enumerate}
\end{lemma}

\begin{proof}
  The only thing that remains to be proved is that $\Psi$ is not
  diagonalizable.  We argue by contradiction.  Suppose there is a
  diagonal function $D$ and fixed unitaries $U,V$ so that $D=U\Psi
  V^*$.  Of course $D$ must be unitary valued on $B$.  In particular,
  $D(1)$ is unitary and so by multiplying on the left (or right) by
  $D(1)^*$, it may be assumed that $D(1)=I$.  Since $\Psi(1)=I$,
  $V=U$.
  
  Let $\varphi_1,\varphi_2$ denote the diagonal entries of $D$.  Since
  $D$ is unitary valued on $B$, both $\varphi_1$ and $\varphi_2$ are
  unimodular on $B$.  Further, as $\det(\Psi)$ has $6$ zeros and a
  function unimodular on $B$ has at least three zeros or is constant,
  we conclude, either both $\varphi_1$ and $\varphi_2$ have three
  zeros and thus take each value in $\{|z|\le 1\}$ exactly three times
  in $X$, or one has six zeros and the other is a unimodular constant
  $\gamma$.  This later case cannot occur, since then
  \begin{equation*}
    0=\Psi(0)=U^* \begin{pmatrix} \gamma & 0\\ 0 & 0 \end{pmatrix}U.
  \end{equation*}
  
  From Lemma \ref{lem:PSIONE}, $\Psi(q_1)e_1=e_1$.  Thus $Ue_1$ is an
  eigenvector for $D(q_1)$ corresponding to eigenvalue $1$, whence at
  least one of $\varphi_j(q_1)$ is $1$.  Similarly, $Ue_2$ is an
  eigenvector for $D(q_1^*)$ with eigenvalue $1$, and at least one of
  $\varphi_j(q_1^*)$ is $1$.  Now, $D(q_1)$ cannot be a multiple of
  the identity, as otherwise both $\varphi_j(q_1)=1$, in which case at
  least one of these two functions takes the value $1$ at both $q_1$
  and $q_1^*$.  Therefore, without loss of generality, $Ue_1=\lambda_1
  e_1$ and $Ue_2=\lambda_2 e_2$, for unimodular $\lambda_1$ and
  $\lambda_2$.  Since $D$ is also diagonal, we can assume
  $\lambda_j=1$.  Hence, $\Psi=D$.
  
  From Lemma \ref{lem:PSIONE}, $\Psi(q_2)(\eta e_1+(1-\eta^2)^\frac12
  e_2)=\eta e_1+(1-\eta^2)^\frac12 e_2$ so that
  \begin{equation*}
    D(q_2)(\eta e_1+(1-\eta^2)^\frac12 e_2)=\eta \varphi_1(q_2)e_1
    +(1-\eta^2)^\frac12 \varphi_2(q_2)e_2 .
  \end{equation*}
  Since $\eta \ne 0$, it follows that $\varphi_j(q_2)=1$.  Similarly,
  $\varphi_j(q_2^*)=1$.  As $q_2\ne q_2^*$, the function $\varphi_j$
  takes the value $1$ twice on $B_2$, a contradiction.
\end{proof}

To prove our main theorem, simply note that the first three conditions
of Lemma \ref{lem:propsoffunnypsi} imply, in view of Theorem
\ref{thm:diagonal}, that if $\rho_\Psi=1$, then $\Psi$ is diagonal.
Hence, it follows that $\rho_\Psi<1$ and this completes the proof.


\section{Existence of a Finite Dimensional Example}


Let $\mathcal A$ denote the closure of $\mathcal R (X)$, the rational
functions with poles off of $X$ as a subspace of $C(X)$.  That is
$\mathcal A=\overline{\mathcal R (X)} \subset C(X)$.  The algebra
$\mathcal A$ is an (abstract) operator algebra with the family of
matrix norms,
\begin{equation*}
  \|F\|=\sup\{\|F(z)\|:z\in X\},
\end{equation*}
for $F=(f_{j,\ell})\in M_n(\mathcal A)$ and where $\|F(z)\|$ is the
matrix norm of the $n\times n$ matrix $F(z)$.
 
Given a finite subset $\Lambda \subset R$, let $I_\Lambda$ denote the
ideal $\{f\in A:f(\lambda)=0 \text{ for all } \lambda\in\Lambda\}$ of
$\mathcal A$.  The quotient
\begin{equation*}
  \mathcal A_\Lambda = \mathcal A/ I_\Lambda,
\end{equation*} 
inherits, in a canonical way, an operator algebra structure from
$\mathcal A$.  Namely,
\begin{equation*}
  \|\pi_\Lambda(F)\|=\inf\left\{ \|G\|: G\in M_n(\mathcal A) \text{
      and } \pi_\Lambda(G)=\pi_\Lambda(F)\right\} .
\end{equation*}
Here $\pi_\Lambda$ denote the quotient map.

The following is a special case of a version of a theorem of Vern
Paulsen \cite{VIP}.  A representation $\nu$ of $\mathcal A$ is a
unital homomorphism $\nu:\mathcal A \to B(H)$, where $H$ is a Hilbert
space.  The representation $\nu$ is contractive if $\|\nu(f)\|\le
\|f\|$ for all $f\in F$ and is completely contractive if for each $n$
and each $F=(f_{j,\ell})\in M_n(\mathcal A)$, $\|\nu(F)\|=
\|(\nu(f_{j,\ell}))\|\le \|F\|$.  Make similar definitions for
representations, contractive representations, and completely
contractive representations of $\mathcal A_E$.

\begin{theorem}[Paulsen]
  Every contractive representation of $\mathcal A$ is completely
  contractive if and only if every contractive representation of every
  $\mathcal A_\Lambda$ is completely contractive.
\end{theorem}

The following observation of Paulsen is also useful.  We give the
proof.

\begin{proposition}
  If $\nu:\mathcal A_\Lambda \to B(H)$ is contractive, but not
  completely contractive, then $\nu\circ\pi_\Lambda:\mathcal A\to
  B(H)$ is contractive, but not completely contractive.
\end{proposition}

\begin{proof}
  Since $\nu$ and $\pi_\Lambda$ are contractive, the composition
  $\nu\circ\pi_\Lambda$ is contractive.  On the other hand, if
  $\nu\circ\pi_\Lambda$ is completely contractive, then it follows
  from the definition of the norm on $M_n(\mathcal A_E)$ as the
  quotient norm, that $\nu$ is completely contractive.  To see this,
  given $\pi_\Lambda(F)=(\pi_\Lambda (f_{j,\ell}))$ is in
  $M_n(\mathcal A_E)$ we can assume, by the definition of the quotient
  norm, that $\|F\|$ is only a little larger than
  $\|\pi_\Lambda(F)\|$.  Since $\nu\circ\pi_\Lambda$ is completely
  contractive, the norm of $\nu\circ\pi_\Lambda(F)$ is no larger than
  $\|F\|$ and the result follows.
\end{proof}

We use Paulsen's results, together with our own, to argue that there
is a matrix $X$ with $R$ as a spectral set, but which does not dilate
to a normal operator with spectrum in $B$.  It is enough to prove that
there is a finite subset $\Lambda\subset R$, a finite dimensional
Hilbert space $H$, and a representation $\nu:\mathcal A_E\to B(H)$
which is contractive, but not completely contractive.  We begin with a
lemma.

\begin{lemma}
  Fix $\Lambda$ a finite subset of $R$.  If $\gamma:\mathcal
  A_\Lambda\to B(K)$ is a representation, then there exists subspaces
  $\mathcal E_\lambda$, $\lambda\in\Lambda$, of $K$ so that
   \begin{enumerate}
   \item each $\mathcal E_\lambda$ is closed;
   \item $K$ is the algebraic direct sum $\hat{\oplus} \mathcal
     E_\lambda$; and
   \item if $f\in\mathcal A$ and $k\in \mathcal \mathcal E_\lambda$,
     then
     \begin{equation*}
       \gamma(\pi_\Lambda(f))k=\lambda k.
     \end{equation*}
   \end{enumerate}
\end{lemma}

\begin{proof}
  For each $\lambda\in \Lambda$, choose a function
  $e_\lambda\in\mathcal A$ such that $e_\lambda(\mu)$ is $0$ if
  $\mu\in \Lambda$, but $\mu\ne \lambda$ and such that
  $e_\lambda(\lambda)=1$.  It follows that the operators $E_\lambda
  =\gamma(\pi_\Lambda(e_\lambda))$ are idempotents such that
  $E_\lambda E_\mu=0$ if $\mu\ne \lambda$ and $\sum E_\lambda$ is the
  identity on $K$.  Let $\mathcal E_\lambda =E_\lambda K$.  Properties
  (1) and (2) are readily verified.  As for (3),
  \begin{equation*}
    \begin{split}
      (\gamma(\pi_\Lambda (f))-f(\lambda))E_\lambda
      &= \gamma(\pi_\Lambda (f-f(\lambda))e_\lambda)\\
      &= \gamma(\pi_\Lambda (0))=0.
    \end{split}
  \end{equation*}
\end{proof}

From our result, there is a representation of $\mathcal A$ which is
contractive, but not completely contractive.  From Paulsen's theorem,
there exists a $\Lambda$ and Hilbert space $K$ (possibly infinite
dimensional) and a representation $\gamma:\mathcal A_\Lambda \to B(K)$
which is contractive, but completely contractive.  Since this
representation is not completely contractive, there exists an $n$ and
$F\in M_n(\mathcal A)$ such that
\begin{equation*}
  1=\|\pi_\Lambda(F)\|<\|\gamma(\pi_\Lambda(F))\|.
\end{equation*}
Thus, there exists a vector
\begin{equation*}
  x=\begin{pmatrix} x_1\\x_2\\ \vdots \\ x_n \end{pmatrix} \in
  \oplus_1^n K=\mathbb C^n\otimes K
\end{equation*}
such that $\|x\|=1$ and
\begin{equation*}
  \| \gamma(\pi_\Lambda(F)) x\|>1.
\end{equation*} 
Write, $x_j =\sum x_j(\lambda)$ with respect to the algebraic direct
sum $\hat{\oplus}\mathcal E_\lambda$ and let $\mathcal H$ denote the
span of $\{x_j(\lambda):1\le j\le n, \ \ \lambda\in \Lambda\}$.  Thus,
$\mathcal H$ is finite dimensional (the dimension is cardinality of
$\Lambda$ times $n$).  From condition (3) of our lemma, $\mathcal H$
is invariant for $\gamma$.  Define $\nu:\mathcal A_E\to B(\mathcal H)$
by restriction: $\nu(f)=\gamma(f)|_{\mathcal H}$.  Since also $x\in
\mathbb C^n\otimes \mathcal H$ we have, $\gamma(\pi_\Lambda(F))
x=\nu(\pi_\Lambda(F)) x$ and thus
\begin{equation*}
  \| \nu(\pi_\Lambda(F)) x\|>1= \|\pi_\Lambda(F)\|.
\end{equation*} 

\bibliographystyle{plain} \bibliography{rat_dil}

\end{document}